# NONPARAMETRIC ESTIMATION OF COMPOSITE FUNCTIONS[1]

By Anatoli B. Juditsky, Oleg V. Lepski
and Alexandre B. Tsybakov

*Université Grenoble 1, Université de Provence and
CREST and Université Paris 6*

We study the problem of nonparametric estimation of a multivariate function $g:\mathbb{R}^d \to \mathbb{R}$ that can be represented as a composition of two unknown smooth functions $f:\mathbb{R}\to\mathbb{R}$ and $G:\mathbb{R}^d \to \mathbb{R}$. We suppose that $f$ and $G$ belong to known smoothness classes of functions, with smoothness $\gamma$ and $\beta$, respectively. We obtain the full description of minimax rates of estimation of $g$ in terms of $\gamma$ and $\beta$, and propose rate-optimal estimators for the sup-norm loss. For the construction of such estimators, we first prove an approximation result for composite functions that may have an independent interest, and then a result on adaptation to the local structure. Interestingly, the construction of rate-optimal estimators for composite functions (with given, fixed smoothness) needs adaptation, but not in the traditional sense: it is now adaptation to the local structure. We prove that composition models generate only two types of local structures: the local single-index model and the local model with roughness isolated to a single dimension (i.e., a model containing elements of both additive and single-index structure). We also find the zones of $(\gamma, \beta)$ where no local structure is generated, as well as the zones where the composition modeling leads to faster rates, as compared to the classical nonparametric rates that depend only to the overall smoothness of $g$.

**1. Introduction.** In this paper we study the problem of nonparametric estimation of an unknown function $g:\mathbb{R}^d \to \mathbb{R}$ in the multidimensional Gaussian white noise model described by the stochastic differential equation

(1) $\qquad X_\varepsilon(dt) = g(t)\,dt + \varepsilon W(dt), \qquad t = (t_1,\ldots,t_d) \in \mathcal{D},$

where $\mathcal{D}$ is a bounded open interval in $\mathbb{R}^d$ containing $[-1,1]^d$, $W$ is the standard Brownian sheet in $\mathbb{R}^d$ and $0 < \varepsilon < 1$ is a known noise level. Our goal is to estimate the function $g$ on the set $[-1,1]^d$ from the observation $\{X_\varepsilon(t), t \in \mathcal{D}\}$. For $d = 2$ this corresponds to the problem of image reconstruction from observations corrupted by additive noise. We consider observation set $\mathcal{D}$, which is larger than $[-1,1]^d$ in order to avoid the discussion of boundary effects.

To measure the performance of estimators, we use the risk function determined by the sup-norm $\|\cdot\|_\infty$ on $[-1,1]^d$: for $g:\mathbb{R}^d \to \mathbb{R}$, $0 < \varepsilon < 1$, $p > 0$, and for an arbitrary estimator $\tilde{g}_\varepsilon$ based on the observation $\{X_\varepsilon(t), t \in \mathcal{D}\}$ we consider the risk

$$(2) \qquad R_\varepsilon(\tilde{g}_\varepsilon, g) = \mathbb{E}_g(\|\tilde{g}_\varepsilon - g\|_\infty^p).$$

Here and in what follows $\mathbb{E}_g$ denotes the expectation with respect to the distribution $\mathbb{P}_g$ of the observation $\{X_\varepsilon(t), t \in \mathcal{D}\}$ satisfying (1).

We suppose the $g \in \mathcal{G}_\mathsf{s}$, where $\{\mathcal{G}_\mathsf{s}, \mathsf{s} \in \mathsf{S}\}$ is a collection of functional classes indexed by $\mathsf{s} \in \mathsf{S}$. The functional classes $\mathcal{G}_\mathsf{s}$ that we will consider consist of *smooth composite functions* and below we discuss in detail this choice.

For a given class $\mathcal{G}_\mathsf{s}$ we define the maximal risk

$$(3) \qquad R_\varepsilon(\tilde{g}_\varepsilon, \mathcal{G}_\mathsf{s}) = \sup_{g \in \mathcal{G}_\mathsf{s}} R_\varepsilon(\tilde{g}_\varepsilon, g).$$

Our first aim is to study the asymptotics, as the noise level $\varepsilon$ tends to 0, of the minimax risk

$$\inf_{\tilde{g}_\varepsilon} R_\varepsilon(\tilde{g}_\varepsilon, \mathcal{G}_\mathsf{s}),$$

where $\inf_{\tilde{g}_\varepsilon}$ denotes the infimum over all estimators of $g$. We suppose that parameter $\mathsf{s}$ is known, and therefore the functional class $\mathcal{G}_\mathsf{s}$ is fixed. We find the minimax rate of convergence $\phi_\varepsilon(\mathsf{s})$ on $\mathcal{G}_\mathsf{s}$, that is, the rate that satisfies $\phi_\varepsilon^p(\mathsf{s}) \asymp \inf_{\tilde{g}_\varepsilon} R_\varepsilon(\tilde{g}_\varepsilon, \mathcal{G}_\mathsf{s})$, and we construct an estimator attaining this rate, which we refer to as a rate-optimal estimator in the asymptotic minimax sense.

**2. Global rate-optimal estimation via pointwise selection.** In this section we discuss a rather general method of data-driven selection from a given family of estimators. This method, called a *pointwise selection rule*,[2] is at the core of the paper. We will use it to construct our rate-optimal estimators.

---

[2] This selection rule was the topic of the IMS Medallion Lecture given by the second author at the *Joint Statistical Meetings* in Minneapolis, 2005.

To present the pointwise selection rule we need some definitions. Let $\mathcal{D}_1$ be an open interval such that $[-1,1]^d \subset \mathcal{D}_1 \subset \mathcal{D}$. Any function $K:\mathbb{R}^d \times \mathbb{R}^d \to \mathbb{R}$ such that

$$\int_{\mathcal{D}_1} K(t,x)\,dt = 1 \qquad \forall x \in [-1,1]^d,$$

$$\operatorname{supp} K(\cdot, x) \subseteq \mathcal{D} \qquad \forall x \in \mathcal{D}_1,$$

will be called a *weight*. Let $\mathcal{K}$ be a given family of weights and let $x \in [-1,1]^d$ be fixed. To any $K \in \mathcal{K}$ we associate a *linear estimator* at $x$:

$$\hat{g}_K(x) = \int_{\mathcal{D}} K(t,x) X_\varepsilon(dt).$$

We consider a family of linear estimators $\mathcal{G}(\mathcal{K}) = \{\hat{g}_K(x), K \in \mathcal{K}\}$. Note that $\hat{g}_K(x)$ is a normal random variable with variance $\varepsilon^2 \|K(\cdot, x)\|_2^2$ where $\|\cdot\|_2$ denotes the $L_2$ norm. Define $\sigma_K = \sup_{x \in \mathcal{D}} \|K(\cdot, x)\|_2$ and assume that the family $\mathcal{K}$ satisfies:

$$\sup_{K \in \mathcal{K}} \sigma_K < \infty.$$

For any pair of weights $K_1$ and $K_2$ define the function

$$[K_1 \otimes K_2](\cdot, \cdot) = \int_{\mathcal{D}_1} K_1(\cdot, y) K_2(y, \cdot)\, dy.$$

We say that $\mathcal{K}$ is a *commutative weight system* if

$$[K_1 \otimes K_2] = [K_2 \otimes K_1] \qquad \forall K_1,\ K_2 \in \mathcal{K}.$$

We now present the pointwise selection rule and briefly discuss some examples where it can be applied. The rule consists of the following two steps:

1. *Determination of acceptable weights.* Let $\mathcal{K}$ be a commutative weight system and let $\mathbf{th}_\varepsilon(\mathcal{K})$ be a threshold whose choice will be discussed below. We say that a weight $K \in \mathcal{K}$ [resp., the estimator $\hat{g}_K(x)$] is *acceptable* if

   $$|\hat{g}_{K \otimes \tilde{K}}(x) - \hat{g}_{\tilde{K}}(x)| \leq M(\mathcal{K}) \mathbf{th}_\varepsilon(\mathcal{K}) \sigma_{\tilde{K}} \qquad \forall \tilde{K} \in \mathcal{K}: \sigma_{\tilde{K}} \geq \sigma_K,$$

   where $M(\mathcal{K}) = \sup_{K \in \mathcal{K}} \sup_{x \in \mathcal{D}} \|K(\cdot, x)\|_1$ and $\|\cdot\|_1$ denotes the $L_1$ norm.
2. *Selection from the set of acceptable estimators.* Let $\hat{\mathcal{K}}$ be the set of all the acceptable weights in $\mathcal{K}$. Note that $\hat{\mathcal{K}}$ is a random set and it can be empty with some probability. If $\hat{\mathcal{K}} \neq \varnothing$ we select the estimator $\hat{g}_{\hat{K}}(x)$ with $\hat{K}$ such that $\sigma_{\hat{K}} = \inf_{K \in \hat{\mathcal{K}}} \sigma_K$, that is, we choose an acceptable estimator with minimal variance. If $\hat{\mathcal{K}} = \varnothing$ we select an arbitrary fixed estimator $\hat{g}_{K_0}(x)$, where $K_0$ is a given weight from $\mathcal{K}$.

There is no general receipt for the choice of the threshold $\mathbf{th}_\varepsilon(\mathcal{K})$. It may depend on the weight system, on the nature of the considered problem (pointwise or global estimation), on the loss functional, etc. However, if we consider the risk (2) and if the weight system $\mathcal{K}$ is not too large (e.g., $\mathcal{K}$ is a metric compact with a polynomial behavior of covering numbers) it can be shown that there is a *universal* choice of the threshold: $\mathbf{th}_\varepsilon(\mathcal{K}) = C\varepsilon\sqrt{\ln 1/\varepsilon}$, where $C > 0$ is a constant depending only on the power $p$ of the loss function and on the dimension $d$. Such a choice of the threshold will be used in this paper.

A remarkable property of the pointwise selection rule is that it can be shown to work for any commutative weight system. As we will see in the following examples, the commutativity property is inherent to a variety of weight systems used in statistics.

*Examples of commutative weight systems.* We now consider some examples of commutative weight systems. Let $\mathcal{Q}$ be any set of functions $Q:\mathbb{R}^d \to \mathbb{R}$ such that $\mathrm{supp}(Q) \subset [-\delta,\delta]^d, \delta > 0$, and $\int_{\mathbb{R}^d} Q = 1$. Take $\mathcal{D} = [-a,a]^d$ and $\mathcal{D}_1 = [-b,b]^d$, where $a > b > 1, a - b > \delta$ are given numbers. Define

$$\mathcal{K} = \{K:\mathbb{R}^d \times \mathbb{R}^d \to \mathbb{R} : K(t,x) = Q(t-x), Q \in \mathcal{Q}\}.$$

Then $\mathcal{K}$ is a commutative weight system. Indeed, the integration over $\mathcal{D}_1$ in the definition of the weight and in the definition of $[K_1 \otimes K_2]$ can be replaced by integration over $\mathbb{R}^d$, and the operation $\otimes$ reduces to the standard convolution:

$$[K_1 \otimes K_2] = K_1 * K_2 = K_2 * K_1 = [K_2 \otimes K_1].$$

This allows us to construct various commutative weight systems. We now consider some of them.

The selection of an estimator from a given family first appeared in the context of *adaptive estimation*. In particular, in [16] a pointwise selection rule was proposed in order to construct pointwise adaptive estimators over a scale of Hölder classes. This method was generalized in [21] to a pointwise selection rule from the collection $\mathcal{G}(\mathcal{K}_{\mathcal{H}_1})$ with the family of weights

$$\mathcal{K}_{\mathcal{H}_1} = \left\{ h^{-1} Q_0\left(\frac{\cdot - x}{h}\right), h \in \mathcal{H}_1 \right\},$$

where $d = 1$, $Q_0 \in \mathcal{Q}$ is a given function, $\mathcal{H}_1 = [h_{\min}, h_{\max}]$ and the numbers $0 < h_{\min} < h_{\max} \leq 1$ are chosen by the statistician. In words, the family $\mathcal{G}(\mathcal{K}_{\mathcal{H}_1})$ consists of kernel estimators with bandwidth varying from $h_{\min}$ to $h_{\max}$. The estimator chosen from the collection $\mathcal{G}(\mathcal{K}_{\mathcal{H}_1})$ in accordance with the pointwise selection rule of [21] is rate optimal over the Besov classes of functions; compare [19].

More recently, pointwise adaptive methods have been developed in dimensions larger than 1. Thus, [14, 15] propose a pointwise selection rule from the collection $\mathcal{G}(\mathcal{K}_{\mathcal{H}_d})$ where

$$\mathcal{K}_{\mathcal{H}_d} = \left\{ \prod_{i=1}^d h_i^{-1} Q_0 \left( \frac{\cdot - x_i}{h_i} \right), (h_1, \ldots, h_d) \in \mathcal{H}_d \right\}.$$

Here the $x_i$ are the components of $x$, and $\mathcal{H}_d = \prod_{i=1}^d [h_{\min}^{(i)}, h_{\max}^{(i)}]$ with the values $0 < h_{\min}^{(i)} < h_{\max}^{(i)} < \infty, i = 1, \ldots, d$, that are chosen by the statistician. The pointwise selection rule of [14] leads to an estimator that is pointwise adaptive over the scale of anisotropic Besov classes [14, 15].

The results of these papers show that pointwise selection is a useful tool for estimation of functions with inhomogeneous smoothness. Another approach to multivariate function estimation is based on structural models. Typical examples are the *single index model* and the *additive model* (see Section 3 for more details). For such models, an important issue is adaptation to the unknown structure, and it can be also carried out via the pointwise selection rule [8]. The weight system used in pointwise selection for the single-index model [8] will also appear in some parts of the present paper. It makes use of the ridge functions. Another system of ridge functions is proposed in [4, 5] for the problem of recovery of functions of two variables with discontinuities along smooth edges and smooth otherwise. Note that the approach of [4, 5] is conceptually different, and does not rely on pointwise selection rules. Examples of more complex commutative weight systems can be found in [8, 20]. Another construction leading to quite an unusual commutative weight system will be given in Section 6.2.

In the present paper we specify the pointwise selection rule for the problem of estimation of composite functions. Our structural assumption is that the function $g : \mathbb{R}^d \to \mathbb{R}$ can be represented as a composition of two unknown smooth functions $f : \mathbb{R} \to \mathbb{R}$ and $G : \mathbb{R}^d \to \mathbb{R}$, that is, $g = f \circ G$.

**3. Why smooth composite functions.** We now discuss why this structural assumption is relevant. We start with the following definition.

DEFINITION 1. Fix $\alpha > 0$ and $L > 0$. Let $\lfloor \alpha \rfloor$ be the largest integer which is strictly less than $\alpha$, and for $\vec{k} = (k_1, \ldots, k_d) \in \mathbb{N}^d$ set $|\vec{k}| = k_1 + \cdots + k_d$. The *isotropic Hölder class* $\mathbb{H}_d(\alpha, L)$ is the set of all functions $G : \mathbb{R}^d \to \mathbb{R}$ having on $\mathbb{R}^d$ all partial derivatives of order $\lfloor \alpha \rfloor$ and such that

$$\sum_{0 \leq |\vec{k}| \leq \lfloor \alpha \rfloor} \sup_{x \in \mathbb{R}^d} \left| \frac{\partial^{|\vec{k}|} G(x)}{\partial x_1^{k_1} \cdots \partial x_d^{k_d}} \right| \leq L,$$

$$(4) \quad \left| G(y) - \sum_{0 \leq |\vec{k}| \leq \lfloor \alpha \rfloor} \frac{\partial^{|\vec{k}|} G(x)}{\partial x_1^{k_1} \cdots \partial x_d^{k_d}} \prod_{j=1}^{d} \frac{(y_j - x_j)^{k_j}}{k_j!} \right| \leq L \|y - x\|^{\alpha}$$
$$\forall x, y \in \mathbb{R}^d,$$

where $x_j$ and $y_j$ are the $j$th components of $x$ and $y$ and $\|\cdot\|$ is the Euclidean norm in $\mathbb{R}^d$.

Parameter $\alpha$ characterizes the isotropic (i.e., the same in each direction) smoothness of function $G$.

Let now $f$ and $G$ be smooth functions such that $f \in \mathbb{H}_1(\gamma, L_1)$ and $G \in \mathbb{H}_d(\beta, L_2)$ where $\gamma, L_1, \beta, L_2$ are positive constants. Here and in what follows $\mathbb{H}_1(\gamma, L_1)$ and $\mathbb{H}_d(\beta, L_2)$ are the Hölder class on $\mathbb{R}$ and the isotropic Hölder class on $\mathbb{R}^d$, respectively. The class of composite functions $g = f(G(x))$ with such $f$ and $G$ will be denoted by $\mathbb{H}(\mathcal{A}, \mathcal{L})$, where $\mathcal{A} = (\gamma, \beta) \in \mathbb{R}_+^2$ and $\mathcal{L} = (L_1, L_2) \in \mathbb{R}_+^2$.

The performance of an estimation procedure will be measured by the sup-norm risk (3) where we set $\mathsf{s} = (\mathcal{A}, \mathcal{L})$ and $\mathcal{G}_\mathsf{s} = \mathbb{H}(\mathcal{A}, \mathcal{L})$.

3.1. *Motivation* I: *models of reduced complexity.* It is well known that the main difficulty in estimation of multivariate functions is the curse of dimensionality: the best attainable rate of convergence of the estimators deteriorates very fast as the dimension grows. To illustrate this effect, suppose, for example, that the underlying function $g$ belongs to $\mathcal{G}_\mathsf{s} = \mathbb{H}_d(\alpha, L), \mathsf{s} = (\alpha, L)$, $\alpha > 0, L > 0$. Then the rate of convergence for the risk (3), uniformly on $\mathbb{H}_d(\alpha, L)$, cannot be asymptotically better than

$$\psi_{\varepsilon, d}(\alpha) = (\varepsilon \sqrt{\ln(1/\varepsilon)})^{2\alpha/(2\alpha+d)}$$

(cf. [6, 12, 13, 23, 25]). This is also the minimax rate on $\mathbb{H}_d(\alpha, L)$; it is attained, for example, by a kernel estimator with properly chosen bandwidth and kernel. More results on asymptotics of the minimax risks in estimation of multivariate functions can be found in [2, 3, 14, 15, 22]. It is clear that if $\alpha$ is fixed and $d$ is large enough this asymptotics is too pessimistic to be used for real data.

At the origin of this phenomenon is the fact that the $d$-dimensional isotropic Hölder class $\mathbb{H}_d(\alpha, L)$ is too massive in terms of its metric entropy. A way to circumvent the curse of dimensionality is to consider models with slimmer functional classes (i.e., classes with smaller metric entropy). There are several ways to do it.

- A first way is to impose a restriction on the smoothness parameter of the functional class. For the class $\mathbb{H}_d(\alpha, L)$, a convenient restriction is to assume that the smoothness $\alpha$ increases with the dimension, and thus the class becomes smaller (its metric entropy decreases). For instance, we can suppose that $\alpha = \kappa d$ with some fixed $\kappa > 0$. Then the dimension disappears from the expression for $\psi_{\varepsilon,d}(\alpha)$, which means that we escape from the curse of dimensionality. However, the condition $\alpha = \kappa d$ or other similar restrictions that link smoothness and dimension are usually difficult to motivate. An interesting related example is given by the class of functions with bounded integrals of the multivariate Fourier transform [1].

- One can also impose a *structural assumption* on the function $g$ to be estimated. Two classical examples are provided by the single-index and additive structures (cf., e.g., [7, 9, 11, 26]).

  The *single-index structure* is defined by the following assumption on $g$: there exist a function $F_0 : \mathbb{R} \to \mathbb{R}$ and a vector $\vartheta \in \mathbb{R}^d$ with $\|\vartheta\| = 1$ such that $g(x) = F_0(\vartheta^T x)$.

  The *additive structure* is defined by the following assumption: there exist functions $F_i : \mathbb{R} \to \mathbb{R}, i = 1, \ldots, d$, such that $g(x) = F_1(x_1) + \cdots + F_d(x_d)$, where $x_j$ is the $j$th component of $x \in \mathbb{R}^d$.

  If we suppose that $F_i \in \mathbb{H}_1(\alpha, L), i = 0, \ldots, d$, then in both cases function $g$ can be estimated with the rate $(\varepsilon\sqrt{\ln(1/\varepsilon)})^{2\alpha/(2\alpha+1)}$, which does not depend on the dimension and coincides with the minimax rate $\psi_{\varepsilon,1}(\alpha)$ of estimation of functions on $\mathbb{R}$.

In general, under structural assumptions the rate of convergence of estimators improves, as compared to the slow $d$-dimensional rate $\psi_{\varepsilon,d}(\alpha)$. For the above examples the rate does not depend on the dimension.

However, it is often quite restrictive to assume that $g$ has some simple structure, such as the single-index or additive one, *on the whole domain of its definition*. In what follows we refer to this assumption as *global structure*.

A more flexible way of modeling is to suppose that $g$ has a *local structure*. For instance, we can assume that $g$ is well approximated by some single-index or additive structure (or by a combination both) in a small neighborhood of a given point $x$. Local structure depends on $x$ and remains unchanged within the neighborhood. Such an approach can be used to model much more complex objects than the global one. However, the form of the $d$-dimensional neighborhood and the local structure should be chosen by the statistician in advance, which makes the local approach rather subjective.

In the present paper we try to find a compromise between the global and local modeling. Our idea is to consider a sufficiently general global model that would generate suitable local structures, and thus would allow us to construct estimators with nice statistical properties. We argue that this program can be realized for global models where the underlying function $g$ is a composition of two smooth functions.

3.2. *Motivation* II: *structure-adaptive estimation.* The problem of estimation of a composite function can be viewed as that of *structural adaptation*. Indeed, let us suppose that the function $G$ is known and $\beta \geq 1$. It is easy to see that in this case the function $g$ can be estimated with the rate $\psi_{\varepsilon,1}(\gamma)$ corresponding to that of estimation of the univariate function $f$ of smoothness $\gamma$.

Thus, the function $G$ can be considered as a functional nuisance parameter characterizing the unknown structure of the function $g$. An important question in this context is: *what is the price to pay for adaptation to the unknown G?*

Note that the composite model is a kind of generalization of the single-index model; instead of the linear function in the latter model we have here a general function $G$. As discussed above, for the single-index model the optimal rate equals to $\psi_{\varepsilon,1}(\gamma)$. We will show that in the general situation when $G$ is nonlinear, the optimal rate of convergence on $\mathbb{H}(\mathcal{A},\mathcal{L})$ [that we denote $\psi_\varepsilon(\mathcal{A})$] is slower than $\psi_{\varepsilon,1}(\gamma)$, that is, $\psi_{\varepsilon,1}(\gamma)/\psi_\varepsilon(\mathcal{A}) \to 0, \varepsilon \to 0$.

It is easy to see that the class $\mathbb{H}(\mathcal{A},\mathcal{L})$ is contained in the Hölder class $\mathbb{H}_d(\alpha_{\gamma,\beta}, L_3)$, where $L_3 = L_3(\mathcal{L})$ and

$$\alpha_{\gamma,\beta} \triangleq \begin{cases} \gamma\beta, & \text{if } 0 < \gamma, \beta \leq 1, \\ \min(\gamma,\beta), & \text{otherwise.} \end{cases}$$

This inclusion implies that if we ignore the composition structure, that is, if we simply suppose that $g \in \mathbb{H}(\alpha_{\gamma,\beta}, L_3)$, then we can only guarantee the rate of convergence $\psi_{\varepsilon,d}(\alpha_{\gamma,\beta})$. On the other hand, it follows from our results given below that $\psi_\varepsilon(\mathcal{A})/\psi_{\varepsilon,d}(\alpha_{\gamma,\beta}) \to 0, \varepsilon \to 0$, for various values of the regularity parameter $\mathcal{A}$. In other words, the knowledge of the fact that we have a composition structure allows us to improve the rate of convergence as compared to the rate of the best estimator, which only relies on the smoothness properties of $g$.

However, for certain values of the parameter $\mathcal{A} = (\gamma, \beta)$ no improvement due to the structure can be expected. This happens when the structural assumption is essentially equivalent to the fact that $g$ belongs to some isotropic Hölder class. This effect takes place for the following values of $(\gamma, \beta) \in \mathbb{R}^2$:

1°. $0 < \gamma, \beta \leq 1$ *(zone of slow rate).* Clearly, in this zone $\mathbb{H}(\mathcal{A},\mathcal{L}) \subset \mathbb{H}_d(\gamma\beta, L_3)$, where $L_3$ is a positive constant depending only on $\gamma, \beta$ and $\mathcal{L}$. Due to this inclusion a standard kernel estimator with properly chosen bandwidth and the boxcar kernel converges with the rate $\psi_{\varepsilon,d}(\gamma\beta) = (\varepsilon\sqrt{\ln(1/\varepsilon)})^{2\gamma\beta/(2\gamma\beta+d)}$. It is not hard to see (cf. Theorem 1) that this rate is optimal, that is, that a lower bound on the minimax risk holds with the same "slow" rate $\psi_{\varepsilon,d}(\gamma\beta)$ (note that $\gamma\beta \leq 1$).

$2°$. $\gamma \geq \beta, \gamma \geq 1$ *(zone of inactive structure).* In this zone we easily get the inclusions $\mathbb{H}_d(\beta, L_4) \subset \mathbb{H}(\mathcal{A}, \mathcal{L}) \subset \mathbb{H}_d(\beta, L_5)$, where $L_4$ and $L_5$ are positive constants depending only on $\beta$ and $\mathcal{L}$. To show the left inclusion it suffices to consider a set of composite functions with linear $f$ and $G \in \mathbb{H}_d(\beta, L)$. Therefore, the asymptotics of the minimax risk on $\mathbb{H}(\mathcal{A}, \mathcal{L})$ is the same as for an isotropic Hölder class $\mathbb{H}_d(\beta, \cdot)$, that is, the minimax rate on this class is $\psi_{\varepsilon,d}(\beta)$. Note that here we estimate as if there were no structure, and the asymptotics of the minimax risk does not depend on $\gamma$. This explains why we refer to this zone as that of *inactive structure*.

We finally remark that if $\beta \leq 1$ the composite function $g$ is rather non-smooth. The effective smoothness equals to $(1 \wedge \gamma)\beta$, and in view of the above discussion, the minimax rate of convergence of estimators on $\mathbb{H}(\mathcal{A}, \mathcal{L})$ is the same as on the Hölder class $\mathbb{H}_d((1 \wedge \gamma)\beta, \cdot)$. This is a very slow rate $\psi_{\varepsilon,d}((1 \wedge \gamma)\beta)$. Therefore, only for $\beta > 1$ one can expect to find estimators with interesting statistical properties.

**4. Main results.** In this section we state the main results and outline the estimation method. The formal description of the estimation procedure and the proofs are deferred to Sections 5 and 7.1–7.2, respectively.

4.1. *Lower bound for the risks of arbitrary estimators.* For any $\mathcal{A} = (\gamma, \beta) \in \mathbb{R}_+^2$ define

(5) $$\phi_\varepsilon(\gamma, \beta) = \begin{cases} (\varepsilon\sqrt{\ln(1/\varepsilon)})^{2\gamma/(2\gamma+1+(d-1)/\beta)}, \\ \quad \text{if } \beta > 1, \beta \geq d(\gamma - 1) + 1, \\ (\varepsilon\sqrt{(\ln 1/\varepsilon)})^{2/(2+d/\beta)}, \\ \quad \text{if } \gamma > 1, \beta < d(\gamma - 1) + 1, \\ (\varepsilon\sqrt{\ln(1/\varepsilon)})^{2/(2+d/(\gamma\beta))}, \\ \quad \text{if } (\gamma, \beta) \in (0, 1]^2. \end{cases}$$

The boundaries between the zones of these three different rates in $\mathbb{R}_+^2$ are presented by the dashed lines in Figure 1.

An asymptotic lower bound for the minimax risk on $\mathbb{H}(\mathcal{A}, \mathcal{L})$ is given by the following theorem.

THEOREM 1. *For any $\mathcal{A} = (\gamma, \beta) \in \mathbb{R}_+^2$ and any $p > 0$ we have*

$$\liminf_{\varepsilon \to 0} \inf_{\tilde{g}_\varepsilon} \sup_{g \in \mathbb{H}(\mathcal{A}, \mathcal{L})} \mathbb{E}_g[(\phi_\varepsilon^{-1}(\gamma, \beta)\|\tilde{g}_\varepsilon - g\|_\infty)^p] > 0,$$

*where $\inf_{\tilde{g}_\varepsilon}$ denotes the infimum over all estimators of $g$.*

The theorem states that the rate of convergence $\phi_\varepsilon(\gamma, \beta)$ cannot be improved by any estimator. We will show below that for $0 < \gamma, \beta \leq 2$ there exist estimators attaining this rate. Before proceeding to the corresponding result, we make several remarks on the properties of the rate $\phi_\varepsilon(\gamma, \beta)$.

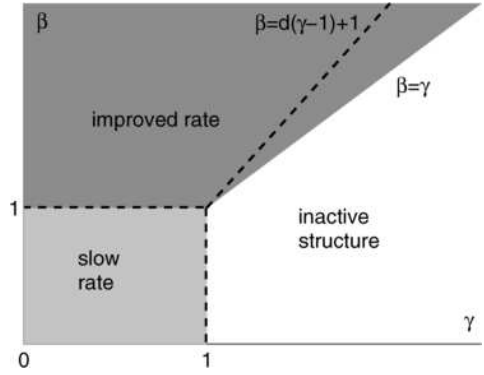

FIG. 1. *Zones of improved rate, of slow rate and of inactive structure. Dashed lines delimit the zones of three different expressions for the rate $\phi_\varepsilon$.*

REMARKS. 1. The set $\{\mathcal{A} = (\gamma, \beta) : \beta > \gamma, \beta \geq 1\}$ will be referred to as the *zone of improved rate* (cf. Figure 1). In this zone there is an improvement of the rate of convergence due to the structure. Indeed, if $\mathcal{A}$ belongs to this zone, the smoothness of function $g$ is equal to $\alpha_{\gamma,\beta} = \gamma$ (cf. Section 3.2), and hence our rate $\phi_\varepsilon(\gamma, \beta)$ is asymptotically (as $\varepsilon \to 0$) much smaller than the rate $\psi_{\varepsilon,d}(\alpha_{\gamma,\beta})$ obtained for the estimators that take into account only the smoothness, and not the structure.

2. The parameter $\beta$ is the tuning parameter of the model: when the ratio $d/\beta$ tends to 0, the rate $\phi_\varepsilon(\gamma, \beta)$, depending on the value of $\gamma$, approaches either the one-dimensional Hölder class rate $\psi_{\varepsilon,1}(\gamma)$ or the "almost parametric" rate $\varepsilon \sqrt{\ln(1/\varepsilon)}$. In particular, when $\beta \geq \gamma > 1$ and $\beta < d(\gamma - 1) + 1$ the rate of convergence $\phi_\varepsilon(\gamma, \beta)$ does not depend on $\gamma$ and coincides with the minimax rate $\psi_{\varepsilon,d}(\beta)$ associated to the $d$-dimensional Hölder class $\mathbb{H}_d(\beta, \cdot)$, and in this zone the composite function $g = f \circ G$ can be estimated with the same rate as $G$, independently of how smooth is $f$.

3. Theorem 1 states the lower bound $(\varepsilon \sqrt{\ln(1/\varepsilon)})^{2\gamma/(2\gamma + 1 + (d-1)/\beta)}$, which is valid for all positive $\gamma, \beta$. Inspection of its proof shows that for $d = 2$ the lower bound is attained on the functions of the form $f_0(\varphi_1(t_1) + \varphi_2(t_2))$. Here $f_0$ is a function of Hölder smoothness $\gamma$ and both functions $\varphi_j, j = 1, 2$, are of Hölder smoothness $\beta$. So, for $d = 2$ the lower bound with the rate $(\varepsilon \sqrt{\ln(1/\varepsilon)})^{2\gamma/(2\gamma + 1 + 1/\beta)}$ holds for that functional family for any $\gamma$ and $\beta$. Note that when $\gamma = \beta$, this lower rate becomes $(\varepsilon \sqrt{\ln(1/\varepsilon)})^{2\beta^2/(2\beta^2 + \beta + 1)}$. Since $\frac{2\beta^2}{2\beta^2 + \beta + 1} < \frac{2\beta}{2\beta + 1}$ this is always slower than the classical one-dimensional rate $\varepsilon^{2\beta/(2\beta + 1)}$. On the other hand, a recent result of [10] shows that for $\gamma = \beta$ functions of the form $f_0(\varphi_1(t_1) + \varphi_2(t_2))$ can be estimated at the rate $\varepsilon^{2\beta/(2\beta + 1)}$ in the $L_2$-norm. Thus, we observe that there is a significant gap between the optimal rates of convergence in $L_2$ and in $L_\infty$, in contrast to the classical nonparametric estimation problems where these rates only differ in a logarithmic factor.

4.2. *Outline of the estimation method.* The exact definition of our estimator is given in Section 5. Here we only outline its construction. We suppose that $\mathcal{A} = (\gamma, \beta) \in (0, 2]^2$. The initial building block is a family of linear estimators. In contrast to the classical kernel construction, which involves a unique bandwidth parameter, the weight $K_\mathcal{J}$ that we consider is determined by the triplet $\mathcal{J} = (\mathcal{A}, \vartheta, \lambda)$ where the *form* parameter $\mathcal{A}$ is the couple $(\gamma, \beta) \in (0, 2]^2$, the *orientation* parameter $\vartheta$ is a unit vector in $\mathbb{R}^d$ and $\lambda$ is a positive real, which we refer to as *size* parameter. We denote $\mathfrak{J}$ the set of all such triplets $\mathcal{J}$ and consider a family of linear estimators $(\hat{g}_\mathcal{J}, \mathcal{J} \in \mathfrak{J})$ where for any $x \in [-1, 1]^d$ the estimator $\hat{g}_\mathcal{J}(x)$ of $g(x)$ is given by

$$\hat{g}_\mathcal{J}(x) \triangleq \int_\mathcal{D} K_\mathcal{J}(t - x) X_\varepsilon(dt).$$

Note that here the size parameter $\lambda$ does not represent the bandwidth of the classical kernel estimator, but rather characterizes the bias of the estimator $\hat{g}_\mathcal{J}$ when the orientation of the window $\vartheta$ is correctly chosen. Namely, the weight $K_\mathcal{J}$ is chosen in such a way that for each $x \in [-1, 1]^d$ the bias of $\hat{g}_\mathcal{J}$ is of the order $O(\lambda)$ if $\vartheta = \vartheta_0^x$ is collinear to the gradient $\nabla G(x)$.

The estimation method proceeds in three steps, and the basic device underlying the construction of the optimal estimation method is the notion of the *local model*. It is an important feature of the composition structure that different local models arise in different subsets of the zone of improved rate.

*Step 1: specifying a collection of local models.* The underlying function $g$ of complicated global structure can have a simple local structure. However, the local structure depends on the function itself. Therefore, $g$ can be only described by a *collection of local models*. In our case, this collection is indexed by a finite-dimensional parameter that can be considered as a nuisance parameter. Specifically, we pass from the global composition model defined in Section 3 to a family of local models $\{\mathcal{M}_\mathcal{J}(x), \mathcal{J} \in \mathfrak{J}, x \in [-1, 1]^d\}$ where the type of each local model $\mathcal{M}_\mathcal{J}(x)$, $\mathcal{J} = (\mathcal{A}, \vartheta, \lambda)$, is determined by $\mathcal{A}$, while $\vartheta$ and $\lambda$ are the local orientation and size parameters. Depending on the value of $\mathcal{A} = (\gamma, \beta)$ (cf. Figure 2), our global model induces only two types of local models: a *local single-index model* and the model with *roughness isolated to a single dimension* (local RISD model).

*1°. Local single-index model: $\gamma \leq 1, 1 < \beta \leq 2$.* In this domain of $\gamma, \beta$, using the smoothness properties of functions $f$ and $G$, it is not hard to show that in the ball $B_{\lambda,x}(\mathcal{A}) = \{t \in \mathbb{R}^d : \|t - x\| \leq \lambda^{1/(\gamma\beta)}\}$ the composite function $g(\cdot)$ can be approximated with the accuracy $O(\lambda)$ by the function

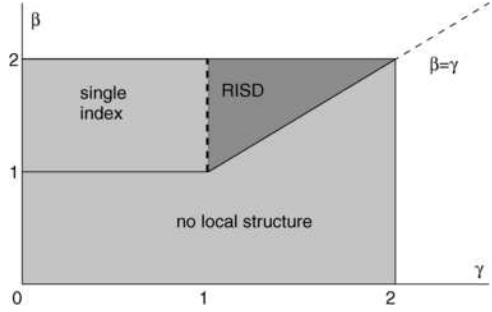

Fig. 2.   *Types of local structures.*

$f(G(x) + \vartheta^T[\cdot - x])$. Here $\vartheta = \vartheta_0^x$ is a unit vector collinear to the gradient $\nabla G(x)$. Indeed, since the inner function $G$ belongs to $H_d(\beta, L_2)$, for any $x, t \in \mathcal{D}$ we have

(6)   $G(t) = G(x) + \nabla G(x)^T(t-x) + B_x(t) \qquad \text{with } |B_x(t)| \leq L_2 \|t-x\|^\beta.$

Next, using the fact that $f \in \mathbb{H}_1(\gamma, L_1)$, we conclude that $g(t) = f(G(t))$ admits the representation

$$g(t) = Q_x(t) + C_x(t),$$

where

$$Q_x(t) = f(G(x) + \nabla G(x)^T(t-x))$$

and

$$|C_x(t)| \leq L_1 |B_x(t)|^\gamma \leq L_1 L_2^\gamma \|t-x\|^{\gamma\beta}.$$

In other words, for any weight $K$ with the support on the ball $B_\lambda(\mathcal{A}) = \{t \in \mathbb{R}^d : \|t\| \leq \lambda^{1/\gamma\beta}\}$ and such that $\int K(y)\,dy = 1$,

(7)   $$\int K(t-x)[g(t) - Q_x(t)]\,dt = O(\lambda).$$

We understand the relation (7) as the definition of the local single-index model $Q_x$ of $g$. The choice of the approximation weight for the function $g$ is naturally suggested by the form of the local model $Q_x$ together with the bound (7): the weight $K_\mathcal{J}$ can be taken as the indicator function of a hyperrectangle normalized by its volume and oriented in such a way that $\nabla G(x)$ is collinear to the first basis vector in $\mathbb{R}^d$. The sides of the hyperrectangle are chosen to have the lengths $l_1 = \lambda^{1/\gamma}$ and $l_j = \lambda^{1/(\gamma\beta)}, j = 2, \ldots, d-1$.

$2°$. *Local model with roughness isolated to a single dimension (RISD):* $1 < \gamma \leq \beta \leq 2$. Let $M_\vartheta$ be an orthogonal matrix with the first column equal to $\vartheta = \vartheta_0^x$, and let $y = M_\vartheta^T(t - x), t \in \mathbb{R}^d$. We denote $y_j$ the $j$th component of $y$ and consider the set

(8) $\quad \mathcal{X}_{\lambda,x}(\mathcal{A}) = \{t \in \mathbb{R}^d : |y_1| \leq \lambda^{1/\beta},\ \|y\| \leq \lambda^{1/(\gamma\beta)},\ |y_1|^{\gamma-1}\|y\|^\beta \leq \lambda\}.$

We show that the estimation of the composite function $g$ at $x$ can be reduced to the problem of estimation under the local model

$$Q_x(y) = q_x(y_1) + P_x(y_2, \ldots, y_d),$$

where $q_x \in \mathbb{H}_1(\gamma, L_1 L_2^\gamma)$ and $P_x \in \mathbb{H}_{d-1}(\beta, 2L_1 L_2)$ on the set $\mathcal{X}_{\lambda,x}(\mathcal{A})$. This local model is established in an unknown coordinate system determined by the parameter $\vartheta = \vartheta_0^x$. Since the smoothness $\gamma$ of $q_x$ is smaller than the smoothness $\beta$ of $P_x$, the accuracy of estimation that corresponds to the coordinate $y_1$ is coarser than that for other coordinates. This motivates the name *roughness isolated to a single dimension*.

The explanation of the local model represented by $Q_x$ on the set $\mathcal{X}_{\lambda,x}(\mathcal{A})$ is provided by the following argument. Using the smoothness properties of functions $f$ and $G$, we obtain due to the inclusions $f \in \mathbb{H}_1(\gamma, L_1)$, $G \in \mathbb{H}_d(\beta, L_2)$:

$$g(t) = f(G(x) + \nabla G(x)^T(t-x)) + f'(G(x) + \nabla G(x)^T(t-x))B_x(t) + C_x(t)$$
$$= f(G(x) + \nabla G(x)^T(t-x)) + f'(G(x))B_x(t) + D_x(t) + C_x(t),$$

where

$$|C_x(t)| \leq C(L_1, L_2, \gamma)\|t - x\|^{\gamma\beta},$$

$$|D_x(t)| \leq C(L_1, L_2)\frac{|\nabla G(x)^T(t-x)|^{\gamma-1}}{\|\nabla G(x)\|}\|t - x\|^\beta, \quad \text{if } \nabla G(x) \neq 0$$

[we have $D_x(t) = 0$ when $\nabla G(x) = 0$], and the function $B_x(t)$, which is defined in (6), belongs to the class $\mathbb{H}_d(\beta, 2L_2)$. In the transformed coordinates (determined by the orthogonal matrix $M_\vartheta$) we may write

(9) $\quad g(t) = g(x + M_\vartheta y) = q(y_1) + \tilde{B}_x(y) + \tilde{D}_x(y) + \tilde{C}_x(y),$

where

(10) $\quad |\tilde{D}_x(y) + \tilde{C}_x(y)| \leq C(L_1, L_2, \gamma)(|y_1|^{\gamma-1}\|y\|^\beta + \|y\|^{\gamma\beta})$

and $\tilde{B}_x \in \mathbb{H}_d(\beta, 2L_2)$. The latter inclusion leads to

(11) $\quad \left|\tilde{B}_x(y) - P_x(y_2, \ldots, y_d) - y_1\frac{\partial}{\partial y_1}\tilde{B}_x(0, y_2, \ldots, y_d)\right| \leq 2L_2|y_1|^\beta,$

where $P_x(y_2, \ldots, y_d) = \tilde{B}_x(0, y_2, \ldots, y_d)$. Let again $K$ be a weight such that $\int K(t)\, dt = 1$, supported on $\mathcal{X}_{\lambda,x}(\mathcal{A})$. Then

$$(12) \qquad \int K(y-x)[g(x + M_\vartheta y) - Q_x(y)]\, dy = O(\lambda),$$

if $K$ is *symmetric in* $y_1$. We understand this property as the definition of the RISD local model $Q_x$ for the composite function $g$.

We conclude that if $\mathcal{A}$ belongs to the zone marked as "RISD" in Figure 2, the *global* structural assumption that the underlying function is a composite one leads automatically to a *local* RISD structure.

A good weight $K_{\mathcal{J}}$ for the zone of RISD local model should be supported on the right window $\mathcal{X}_{\lambda,x}(\mathcal{A})$, possess small bias on both single-index component $q_x$ and "regular" component $P_x$ and have a small $L_2$-norm to ensure small variance of the stochastic term of the estimation error. The construction of such a weight is rather involved (cf. Section 6.2). Note that using a rectangular weight, as for the local single-index model leads to suboptimal estimation rates.

As we see, the definition of local model has two ingredients: the neighborhood (window) and the local structure within the window. For the local single-index model the window is just an Euclidean ball, whereas for the RISD local model the window is the set $\mathcal{X}_{\lambda,x}(\mathcal{A})$.

*Step 2: optimizing the size parameter and specifying candidate estimators.* Once the local model is determined and the corresponding weight is constructed we can choose the size parameter $\lambda = \lambda_\varepsilon(\mathcal{A})$ in an optimal way. To do it we optimize our sup-norm risk with respect to $\lambda$, that is, we get the value $\lambda$, which realizes the balance of bias and variance terms of the risk in the ideal case where the orientation $\vartheta = \vartheta_0^x$ is "correct" for all $x$.

Recall that the weight $K_{\mathcal{J}}$ supported on the window is chosen in such a way that the bias of the linear estimator $\hat{g}_{\mathcal{J}}$, for the "correct" orientation $\vartheta$, is of the order $O(\lambda)$ on every local model. Thus, the bias-variance balance relation for the sup-norm loss can be written in the form

$$(13) \qquad \lambda \asymp \varepsilon \sqrt{\ln 1/\varepsilon} \|K_{\mathcal{J}}\|_2.$$

We will see that $\|K_{\mathcal{J}}\|_2$ depends on $\mathcal{A}$ and $\lambda$ but does not depend on $\vartheta$. This will allow us to choose the optimal value $\lambda_\varepsilon(\mathcal{A})$ independent of $\vartheta$. For instance, for the local single-index model (when $\gamma \leq 1$) the weight $K_{\mathcal{J}}$ is just a properly scaled and rotated indicator of a hyperrectangle. In this particular case the bias-variance balance (13) can be written in the form

$$\lambda \asymp \frac{\varepsilon\sqrt{\ln 1/\varepsilon}}{\sqrt{\text{volume of hyperrectangle}}} = \varepsilon \left( \frac{\ln 1/\varepsilon}{\lambda^{1/\gamma + (d-1)/(\gamma\beta)}} \right)^{1/2}.$$

Note that in this case $\lambda_\varepsilon(\mathcal{A}) \asymp \phi_\varepsilon(\gamma, \beta)$, where $\phi_\varepsilon(\gamma, \beta)$ is defined in (5).

With $\lambda_\varepsilon(\mathcal{A})$ being chosen, we obtain a family of linear estimators

$$
(14) \qquad \{\hat{g}_{\mathcal{J}}(x), \mathcal{J} = (\mathcal{A}, \vartheta, \lambda_\varepsilon(\mathcal{A})) \in \mathfrak{J}, x \in [-1,1]^d\}.
$$

For a fixed $x \in [-1,1]^d$ this family only depends on two parameters, $\mathcal{A}$ and $\vartheta$.

*Step 3: selection.* We now choose an estimator from the family (14) that corresponds to some $\hat{\mathcal{J}} \in \mathfrak{J}$ selected in a data-dependent way, and define our final estimator as a piecewise-constant approximation of the function $x \mapsto \hat{g}_{\hat{\mathcal{J}}}(x)$. To choose $\hat{\mathcal{J}}$ we apply the pointwise selection procedure presented in Section 2.

We introduce a discrete grid on the unit sphere $\{\vartheta \in \mathbb{R}^d : \|\vartheta\| = 1\}$, and we divide the domain of definition of $x$ into small blocks. For each block, we consider a finite set of estimators $\hat{g}_{\mathcal{J}}(x)$ extracted from the family (14), with $x$, which is fixed as the center $x_0$ of the block and all the $\vartheta$ on the grid. We then select a data-dependent $\hat{\vartheta}$ from the grid applying our aggregation procedure to this finite set. The value of our final estimator $g^*_{\mathcal{A},\varepsilon}$ on this block is constant and is defined as $g^*_{\mathcal{A},\varepsilon}(x) \equiv \hat{g}_{(\mathcal{A},\hat{\vartheta},\lambda_\varepsilon(\mathcal{A}))}(x_0)$. We thus get a piecewise-constant estimator $g^*_{\mathcal{A},\varepsilon}$ on $[-1,1]^d$ that depends only on $\mathcal{A}$ and on the observations (the exact definition of $g^*_{\mathcal{A},\varepsilon}$ is given in Section 5).

REMARKS. In this paper we assume that the smoothness $\mathcal{A} = (\gamma, \beta)$ is known, and we deal only with adaptation to the local structure determined by $\vartheta$. If $\mathcal{A}$ is unknown we need simultaneous adjustment of the estimators to $\mathcal{A}$ and to $\vartheta$, that is, to the smoothness and to the local structure of the underlying function. Note, however, that parameters $\mathcal{A}$ and $\vartheta$ are not independent. In particular, $\mathcal{A}$ determines the form of the neighborhood where we have an *unknown* local structure depending on $\vartheta$. This is important because our construction of the family of estimators $\{\hat{g}_{\mathcal{J}}, \mathcal{J} \in \mathfrak{J}\}$ strongly relies on the local representation of the model. For example, if the family $\{\hat{g}_{\mathcal{J}}, \mathcal{J} \in \mathfrak{J}\}$ does not contain an estimator corresponding to the correct local structure, the choice from this family cannot even guarantee consistency. Another difficulty is that different values of $\mathcal{A}$ can correspond to different *types* of local models (cf. Figure 2). In other words, the problem of adaptive estimation of composite functions turns out to be more involved than the classical adaptation to the unknown smoothness as considered, for example, in [16, 17, 18]. As yet we do not know whether fully adaptive estimation in this context is possible or not.

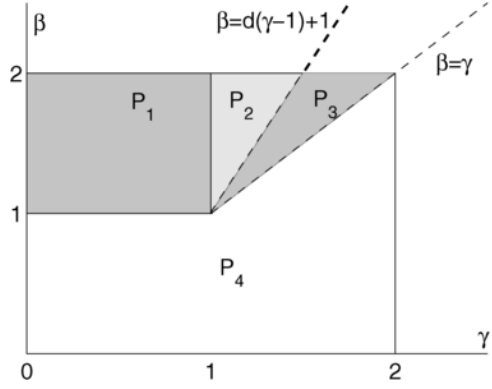

Fig. 3. *Classification of zones within* $(0,2]^2$.

4.3. *Upper bounds on the risk of the estimators.* We define the following three domains of values of $\mathcal{A} = (\gamma, \beta)$ contained in $(0,2]^2$ (cf. Figure 3).

$$
\begin{aligned}
\mathcal{P}_1 &= \{\mathcal{A} : \gamma \leq 1, 1 < \beta \leq 2\}, \\
\mathcal{P}_2 &= \{\mathcal{A} : 1 < \gamma \leq \beta \leq 2, \beta \geq d(\gamma - 1) + 1\}, \\
\mathcal{P}_3 &= \{\mathcal{A} : 1 < \gamma \leq \beta \leq 2, \beta < d(\gamma - 1) + 1\}.
\end{aligned}
\tag{15}
$$

In view of the above discussion, these are exactly the zones where improved rates occur and where the local structure is active. For the sake of completeness, we consider also the remainder zone (zone of no local structure):

$$\mathcal{P}_4 = (0,1]^2 \cup \{(\gamma, \beta) : 1 \leq \beta < \gamma \leq 2\}.$$

As we will see in Section 6.2, the optimal weights $K_{\mathcal{J}}$ are defined separately for each of these zones.

THEOREM 2. *Let $\phi_\varepsilon(\gamma, \beta)$ be as in (5). For any $\mathcal{A} = (\gamma, \beta) \in (0,2]^2 \setminus \mathcal{P}_2$ and any $p > 0$ the estimator $g^*_{\mathcal{A}, \varepsilon}$ satisfies*

$$\limsup_{\varepsilon \to 0} \sup_{g \in \mathbb{H}(\mathcal{A}, \mathcal{L})} \mathbb{E}_g[(\phi_\varepsilon^{-1}(\gamma, \beta) \|g^*_{\mathcal{A}, \varepsilon} - g\|_\infty)^p] < \infty.$$

*For any $\mathcal{A} = (\gamma, \beta) \in \mathcal{P}_2$ and any $p > 0$ the estimator $g^*_{\mathcal{A}, \varepsilon}$ satisfies*

$$\limsup_{\varepsilon \to 0} \sup_{g \in \mathbb{H}(\mathcal{A}, \mathcal{L})} \mathbb{E}_g[([\ln \ln (1/\varepsilon)]^{-1} \phi_\varepsilon^{-1}(\gamma, \beta) \|g^*_{\mathcal{A}, \varepsilon} - g\|_\infty)^p] < \infty.$$

Combining Theorems 1 and 2 we conclude that $\phi_\varepsilon(\gamma, \beta)$ is the minimax rate of convergence for the class $\mathbb{H}(\mathcal{A}, \mathcal{L})$ if $\mathcal{A} = (\gamma, \beta) \in (0,2]^2 \setminus \mathcal{P}_2$, and that it is near minimax [up to the $\ln\ln(1/\varepsilon)$ factor] if $\mathcal{A} = (\gamma, \beta) \in \mathcal{P}_2$. Therefore, our estimator $g^*_{\mathcal{A}, \varepsilon}$ is, respectively, rate optimal or near rate optimal on $\mathbb{H}(\mathcal{A}, \mathcal{L})$.

Theorem 2 is in fact a result on adaptation to the unknown local structure of the function to be estimated: the estimator $g^*_{\mathcal{A},\varepsilon}$ locally adapts to the "correct" orientation $\vartheta_0$, which is collinear to the gradient $\nabla G(x)$ at $x$.

REMARKS. We consider here the Gaussian white noise model because its analysis requires a minimum of technicalities. Composition structures can be studied for more realistic models, such as nonparametric regression with random design, nonparametric density estimation and classification. Note that our theorems can be directly transposed to the Gaussian nonparametric regression model with fixed equidistant design using the equivalence of experiments argument (cf. [24]). Note also that results similar to ours have been recently obtained for the problem of testing hypotheses about composite functions in the Gaussian white noise model [20].

We prove the upper bound of Theorem 2 only for the case $\mathcal{A} \in (0,2]^2$. An extension to $\mathcal{A} \notin (0,2]^2$ remains an open problem. On the other hand, the lower bound of Theorem 1 is valid for all $\mathcal{A} \in \mathbb{R}^2_+$. We believe that it cannot be improved. This conjecture is supported by the recent results on a hypothesis testing problem with composite functions [20], which is closely related to our estimation problem. The upper bound proved in [20] for all $\mathcal{A} \in \mathbb{R}^2_+$ in the problem of hypothesis testing coincides with the lower bound of Theorem 1.

The rate of convergence of the minimax procedure (cf. Theorem 2) in the zone $\mathcal{P}_2$ contains an additional $\ln\ln(1/\varepsilon)$ factor, as compared to the lower bound of Theorem 1. We believe that this minor deterioration of the rate can be avoided by using a more refined estimation procedure.

**5. Definition of the estimator and basic approximation results.** We first introduce some notation. For a bounded function $K \in L_1(\mathbb{R}^d)$ and $p \geq 1$ we denote by $\|K\|_p$ its $L_p$-norm and by $K * g$ its convolution with a bounded function $g$:

$$\|K\|_p = \left(\int |K(t)|^p\, dt\right)^{1/p}, \qquad [K*g](x) = \int K(t-x)g(t)\, dt, \qquad x \in \mathbb{R}^d$$

(here and in the sequel $\int = \int_{\mathbb{R}^d}$). We denote $\mathcal{J} \triangleq (\mathcal{A}, \vartheta, \lambda)$ where $\mathcal{A} = (\gamma, \beta) \in (0,2]^2$, $\vartheta$ is a unit vector in $\mathbb{R}^d$ and $\lambda > 0$. The class of all such triplets $\mathcal{J}$ is denoted by $\mathfrak{J}$.

Given a unit vector $\vartheta$, let $M_\vartheta \in \mathbb{R}^{d \times d}$ stand for an orthogonal matrix with the first column equal to $\vartheta$. The weight system we consider in the sequel is defined as

$$K_\mathcal{J}(x) = \mathsf{K}_{(\mathcal{A},\lambda)}(M_\vartheta^T x),$$

where $\mathsf{K}_{(\mathcal{A},\lambda)}:\mathbb{R}^d\to\mathbb{R}$ is a weight that will be defined in Section 6. Next, for any $\mathcal{J}',\mathcal{J},\in\mathfrak{J}$ and all $t\in\mathbb{R}^d$ we define the *convoluted weight*

$$K_{\mathcal{J}'*\mathcal{J}}(t)=\int K_{\mathcal{J}'}(t-y)K_{\mathcal{J}}(y)\,dy$$

and the difference

$$\Delta_{\mathcal{J}'}K_{\mathcal{J}'*\mathcal{J}}=K_{\mathcal{J}'*\mathcal{J}}-K_{\mathcal{J}'}.$$

We require the weight $K_\mathcal{J}$ to be symmetric, that is, $K_\mathcal{J}(t)=K_\mathcal{J}(-t)$, and

(16) $$K_{\mathcal{J}'*\mathcal{J}}=K_{\mathcal{J}*\mathcal{J}'}.$$

For all $\mathcal{J}\in\mathfrak{J}$ and all $x\in[-1,1]^d$ set

$$\hat{g}_\mathcal{J}(x)=\int_\mathcal{D} K_\mathcal{J}(t-x)X_\varepsilon(dt)$$

and for all $\mathcal{J}',\mathcal{J}\in\mathfrak{J}$ define the *convoluted estimator*

$$\hat{g}_{\mathcal{J}'*\mathcal{J}}(x)=\int_\mathcal{D} K_{\mathcal{J}'*\mathcal{J}}(t-x)X_\varepsilon(dt).$$

In what follows we assume $\varepsilon$ is small enough so that in all expressions that involve weight convolutions we can replace $\int_\mathcal{D}$ by $\int_{\mathbb{R}^d}$ (recall that weights we consider are compactly supported). We also suppose that $\ln\ln(1/\varepsilon)>0$. Define

$$\Delta_{\mathcal{J}'}\hat{g}_{\mathcal{J}'*\mathcal{J}}(x)=\hat{g}_{\mathcal{J}'*\mathcal{J}}(x)-\hat{g}_{\mathcal{J}'}(x)$$

and set

$$\mathbf{TH}_\varepsilon(\mathcal{J}',\mathcal{J})=C(p,d)(\|K_{\mathcal{J}'}\|_1+\|K_\mathcal{J}\|_1)\|K_{\mathcal{J}'}\|_2\varepsilon\sqrt{\ln(1/\varepsilon)},$$

where $C(p,d)=2+\sqrt{4p+8d}$.

5.1. *Estimation procedure.* Now we need to introduce a discrete grid on the set of indices $\mathfrak{J}$. We discretize only the $\vartheta$-coordinate of $\mathcal{J}$. Recall that $\vartheta$ takes values on the Euclidean unit sphere $\mathbb{S}$ in $\mathbb{R}^d$.

*Discretization.* Let $\mathbb{S}_\varepsilon\subset\mathbb{S}$ be an $\varepsilon$-net on $\mathbb{S}$, that is, a finite set such that

$$\forall\vartheta\in\mathbb{S}\qquad\exists\vartheta'\in\mathbb{S}_\varepsilon:\|\vartheta-\vartheta'\|\leq\varepsilon$$

and $\mathrm{card}(\mathbb{S}_\varepsilon)\leq(\sqrt{d}/\varepsilon)^{d-1}$ for small $\varepsilon$. Without loss of generality, we will assume that $(1,0,\ldots,0)\in\mathbb{S}_\varepsilon$.

Fix $\mathcal{A}\in(0,2]^2$ and define $\lambda_\varepsilon(\mathcal{A})$ as a solution in $\lambda$ of the bias-variance balance equation

(17) $$C_1\lambda=\varepsilon\sqrt{\ln(1/\varepsilon)}\|\mathsf{K}_{(\mathcal{A},\lambda)}\|_2,$$

where $C_1$ is a constant in Proposition 2 below, depending only on $\mathcal{A}$, $\mathcal{L}$ and $d$. Finally we define the grid on $\mathcal{J}$:

$$\mathfrak{J}_{\mathrm{grid}}\triangleq\{\mathcal{J}=(\mathcal{A},\vartheta,\lambda_\varepsilon(\mathcal{A})):\vartheta\in\mathbb{S}_\varepsilon\}\subset\mathfrak{J}.$$

*Acceptability.* For a given $x \in [-1,1]^d$ we define a subset $\hat{\mathfrak{T}}_x$ of $\mathfrak{J}_{\text{grid}}$ as follows:

$$\mathcal{J} \in \hat{\mathfrak{T}}_x \iff |\Delta_{\mathcal{J}'}\hat{g}_{\mathcal{J}'*\mathcal{J}}(x)| \leq \mathbf{TH}_\varepsilon(\mathcal{J}', \mathcal{J}) \quad \forall \mathcal{J}' \in \mathfrak{J}_{\text{grid}}.$$

Any value $\mathcal{J} \in \mathfrak{J}_{\text{grid}}$ that belongs to $\hat{\mathfrak{T}}_x$ is called *acceptable*.

Note that the threshold $\mathbf{TH}_\varepsilon(\mathcal{J}', \mathcal{J})$ can be bounded from above and replaced in all the definitions by a value that does not depend on $\mathcal{J}, \mathcal{J}' \in \mathfrak{J}_{\text{grid}}$. In fact, either $\mathbf{TH}_\varepsilon(\mathcal{J}', \mathcal{J}) \asymp \lambda_\varepsilon(\mathcal{A})$ if $\mathcal{A} \in \mathcal{P}_1 \cup \mathcal{P}_3$ or $\mathbf{TH}_\varepsilon(\mathcal{J}', \mathcal{J}) \asymp \ln\ln(1/\varepsilon)\lambda_\varepsilon(\mathcal{A})$ if $\mathcal{A} \in \mathcal{P}_2$.

*Estimation at a fixed point.* For any $x \in [-1,1]^d$ such that $\hat{\mathfrak{T}}_x \neq \varnothing$ we select an arbitrary $\hat{\mathcal{J}}_x$ from the set $\hat{\mathfrak{T}}_x$. Note that the set $\hat{\mathfrak{T}}_x$ is finite, so a measurable choice of $\hat{\mathcal{J}}_x$ is always possible; we assume that such a choice is effectively done. We then define the estimator $g^{**}(x)$ as follows:

$$(18) \quad g^{**}(x) \triangleq \begin{cases} \hat{g}_{\hat{\mathcal{J}}_x}(x), & \text{if } \hat{\mathfrak{T}}_x \neq \varnothing, \\ 0, & \text{if } \hat{\mathfrak{T}}_x = \varnothing. \end{cases}$$

*Global estimator.* The estimator $g^{**}$ is defined for all $x \in [-1,1]^d$ and we could consider $x \mapsto g^{**}(x), x \in [-1,1]^d$, as an estimator of the function $g$. However, the measurability of this mapping is not a straightforward issue. To skip the analysis of measurability, we use again a discretization. Introduce the following cubes in $\mathbb{R}^d$:

$$\Pi_\varepsilon(z) = \bigotimes_{k=1}^d [\varepsilon^2(z_k - 1), \varepsilon^2 z_k], \quad z = (z_1, \ldots, z_d) \in \mathbb{Z}^d.$$

For any $x \in [-1,1]^d$ we consider $z(x) \in \mathbb{Z}^d$ such that $x$ belongs to the cube $\Pi_\varepsilon(z(x))$, and a piecewise constant estimator $g^{**}(z(x))$. Our final estimator is a truncated version of $g^{**}(z(x))$:

$$(19) \quad g^*_{\mathcal{A},\varepsilon}(x) \triangleq \begin{cases} g^{**}(z(x)), & \text{if } |g^{**}(z(x))| \leq \ln\ln(1/\varepsilon), \\ \ln\ln(1/\varepsilon)\operatorname{sign}(g^{**}(z(x))), & \text{if } |g^{**}(z(x))| > \ln\ln(1/\varepsilon). \end{cases}$$

Thus, the resulting procedure $g^*_{\mathcal{A},\varepsilon}$ is piecewise constant on the cubes $\Pi_\varepsilon(z) \subset [-1,1]^d, z \in \mathbb{Z}^d$.

REMARK. Some comments on the numerical complexity of the proposed method are in order here. The algorithm of this section can be easily reformulated for the problem of estimation of the signal $g(i)$ at $n$ points of a regular grid in $[0,1]^d$, from independent observations $y(i) = g(i) + \xi(i)$, $\xi(i) \sim \mathcal{N}(0,1)$, $i = 1, \ldots, n$. A standard argument results in the equivalence between the two models when $\epsilon \asymp n^{-1/2}$, [24].

According to the definition of our method, at each point we need to compare $N = O(n^{(d-1)/2})$ estimators which correspond to the grid over $\vartheta$ on the unit sphere of dimension $d-1$. There are two main components of the numerical effort: we need to compute $N^2$ convoluted weights and the convolutions of these weights with the observation $y$. It will cost $O(n)$ elementary operations to implement the construction of Section 6.2 for each of $N$ weights, and then $O(n \ln n)$ operations to compute each of $N^2$ convolutions. The numerical complexity of this step is therefore $O(N^2 n \ln n) = O(n^d \ln n)$. Further, the convolution of $y$ with each weight requires $O(n \ln n)$ operations. Thus the total cost of convoluting all $N^2$ weights with $y$ will be, again, $O(n^d \ln n)$. Finally, choosing the estimator from the family at each point of the grid demands $N^2$ comparisons. We conclude that the total effort will be $O(n^d \ln n)$ elementary operations, which is far from being prohibitive for dimensions $d=2$ and $d=3$ that are of interest in the context of image analysis.

5.2. *Basic approximation results.* We can now describe the approximation properties of the weight $K_{\mathcal{J}}$, which serve as a main tool in the proof of the properties of the estimator $g^*_{\mathcal{A},\varepsilon}(x)$.

Let $x \in [-1,1]^d$ and $\mathcal{A} = (\gamma, \beta) \in (0,2]^2$ be fixed and let $g = f \circ G \in \mathbb{H}(\mathcal{A}, \mathcal{L})$. We define

$$\vartheta_0^x \triangleq \begin{cases} (1,0,\ldots,0), & \text{if } \beta > 1 \text{ and } \nabla G(x) = 0 \text{ or } \beta \leq 1, \\ \nabla G(x)/\|\nabla G(x)\|, & \text{if } \beta > 1, \nabla G(x) \neq 0. \end{cases} \quad (20)$$

The following statement is an immediate consequence of Lemmas 1–4 formulated in the next section:

PROPOSITION 1. *For all* $\mathcal{A} = (\gamma, \beta) \in (0,2]^2$, *and all* $\lambda > 0$ *we have*

$$\sup_{x \in [-1,1]^d} \sup_{g \in \mathbb{H}(\mathcal{A},\mathcal{L})} |[K_{\mathcal{J}_0^x} * g](x) - g(x)| \leq C_2 \lambda,$$

*where* $\mathcal{J}_0^x = (\mathcal{A}, \vartheta_0^x, \lambda)$ *and* $C_2$ *only depends on* $\mathcal{A}$, $\mathcal{L}$ *and* $d$.

In other words, the weight system $\{K_{\mathcal{J}}, \mathcal{J} \in \mathfrak{J}\}$ contains an element $K_{\mathcal{J}_0^x}$ such that the quality of approximation of $g(x)$ by the "ideal" smoother $[K_{\mathcal{J}_0^x} * g](x)$ is of the order $O(\lambda)$. Here we use the term "ideal" because $\mathcal{J}_0^x = (\mathcal{A}, \vartheta_0^x, \lambda)$ depends on the gradient $\nabla G(x)$, and thus on the unknown function $g$.

The following property of weights $K_{\mathcal{J}}$ is used in the proof of Theorem 2.

PROPOSITION 2. *For all* $\mathcal{A} = (\gamma, \beta) \in (0,2]^2$, $x \in [-1,1]^d$, $0 < \lambda \leq 1$ *and all* $\mathcal{J} = (\mathcal{A}, \vartheta, \lambda) \in \mathfrak{J}$ *we have*

$$\sup_{\mathcal{A} \in (0,2]^2} \sup_{g \in \mathbb{H}(\mathcal{A},\mathcal{L})} |[\Delta_{\mathcal{J}} K_{\mathcal{J}*\mathcal{J}_0^x} * g](x)|$$
$$\leq C_1 \{(\|K_{\mathcal{J}}\|_1 + \|K_{\mathcal{J}_0^x}\|_1)\lambda + \|K_{\mathcal{J}}\|_1 \|K_{\mathcal{J}_0^x}\|_1 \varepsilon\}, \quad (21)$$

where $\mathcal{J}_0^x = (\mathcal{A}, \vartheta^x, \lambda)$, $\vartheta^x$ is any element of the unit sphere $\mathbb{S}$ such that $\|\vartheta^x - \vartheta_0^x\| \leq \varepsilon$ and $C_1$ is a constant depending only on $\mathcal{A}$, $\mathcal{L}$ and $d$. Furthermore, for any $\mathcal{J}, \mathcal{J}' \in \mathfrak{J}$ we have

$$\|\Delta_{\mathcal{J}'} K_{\mathcal{J}' * \mathcal{J}}\|_2 \leq (\|K_{\mathcal{J}'}\|_1 + \|K_{\mathcal{J}}\|_1) \|K_{\mathcal{J}'}\|_2. \tag{22}$$

**6. Weight systems and properties of the weights.** Depending on the value of $\mathcal{A}$ [different zones $\mathcal{P}_i$ (cf. Figure 3)] we use different constructions of $\mathsf{K}_{(\mathcal{A},\lambda)}$. Our objective is to obtain $K_{\mathcal{J}}$ with suitable approximation properties for each $\mathcal{J} \in \mathfrak{J}$. Let us summarize here the main requirements on the weight:

1. Convolution of the weight $\mathsf{K}_{(\mathcal{A},\lambda)}$ with the "local model" of $g$ corresponding to $\mathcal{A}$ should approximate $g$ with the accuracy $O(\lambda)$. Furthermore, the weight should be localized, that is, it should vanish outside of the window where the local structure is valid.
2. A basic characteristic of the weight is its $L_2$-norm, which determines the variance of the estimator. Our objective is to achieve its minimal value.
3. The $L_1$-norm of the weights is also an important parameter of the proposed estimation procedure since it is inherent to the definition of the threshold. Our objective will be to keep the $L_1$-norm as small as possible.

We start with formulation of the properties of the weights, which allows us to prove the basic approximation result and to find the parameters of our estimation procedure. The explicit description of weight systems will be given in the end of the section.

6.1. *Properties of the weights.*

Zone $\mathcal{P}_4$ *(no local structure).*

LEMMA 1. *For any $\mathcal{A} = (\gamma, \beta) \in \mathcal{P}_4$, $\lambda > 0$ and $x \in [-1,1]^d$, we have*

$$\sup_{g \in \mathbb{H}(\mathcal{A},\mathcal{L})} |[\mathsf{K}_{(\mathcal{A},\lambda)} * g](x) - g(x)| \leq c_0 \lambda,$$

*where the constant $c_0$ depends only on $\mathcal{L}$ and $d$. Furthermore,*

$$\|\mathsf{K}_{(\mathcal{A},\lambda)}\|_1 = 1 \quad and \quad \|\mathsf{K}_{(\mathcal{A},\lambda)}\|_2 = \begin{cases} (2\lambda^{1/(\gamma\beta)})^{-d/2}, & (\gamma,\beta) \in (0,1]^2, \\ (2\lambda^{1/\beta})^{-d/2}, & 1 < \beta < \gamma \leq 2. \end{cases}$$

*Zone $\mathcal{P}_1$ (local single-index model).* Let $q : \mathbb{R} \to \mathbb{R}$ and $B : \mathbb{R}^d \to \mathbb{R}$ be functions such that, for given $\gamma \in (0, 1]$,

$$|q(x) - q(y)| \leq L|x - y|^\gamma \qquad \forall x, y \in \mathbb{R}^d,$$

$$\sup_{x \in \mathbb{R}^d} |B(x)| \leq c_1,$$

where $c_1 > 0$, $L > 0$ are constants. We denote by $\mathfrak{A}(\gamma)$ the set of all pairs of functions $(q, B)$ satisfying these restrictions. Define

$$Q(y) = q(y_1) + B(y)\|y\|^{\gamma\beta} \qquad \forall y \in \mathbb{R}^d.$$

We have the following evident result:

LEMMA 2. *For any $\mathcal{A} = (\gamma, \beta) \in \mathcal{P}_1$ and $\lambda > 0$ we have*

(i) $$\sup_{(q,B) \in \mathfrak{A}(\gamma)} |[\mathsf{K}_{(\mathcal{A},\lambda)} * Q](0) - q(0)| \leq c_2 \lambda,$$

*where $c_2$ is a constant depending only on $L$, $c_1$ and $d$. Moreover,*

(ii) $\quad \|\mathsf{K}_{(\mathcal{A},\lambda)}\|_1 = 1 \quad$ and $\quad \|\mathsf{K}_{(\mathcal{A},\lambda)}\|_2 = (2^d \lambda^{1/\gamma + (d-1)/(\gamma\beta)})^{-1/2}.$

*Zone $\mathcal{P}_2 \cup \mathcal{P}_3$ (RISD local model).* Let $q : \mathbb{R} \to \mathbb{R}$ and $p : \mathbb{R}^d \to \mathbb{R}$, $B : \mathbb{R}^d \to \mathbb{R}$ be functions such that $p$ is continuously differentiable and, for given $\mathcal{A} = (\gamma, \beta) \in \mathcal{P}_2 \cup \mathcal{P}_3$ and $\lambda > 0$,

(23) $$\left| q(0) - \frac{1}{2\lambda^{1/\gamma}} \int_{-\lambda^{1/\gamma}}^{\lambda^{1/\gamma}} q(z)\,dz \right| \leq c_3 \lambda,$$

(24) $\quad |p(z') - p(z) - [\nabla p(z)]^T(z' - z)| \leq L\|z' - z\|^\beta \qquad \forall z, z' \in \mathbb{R}^d,$

(25) $$\sup_{x \in \mathbb{R}^d} |B(x)| \leq c_4,$$

where $c_3, c_4$ and $L$ are positive constants. Let $\mathfrak{B}(\mathcal{A}, \lambda)$ denote the set of triplets $(q, p, B)$ satisfying (23)–(25). Define

$$Q(y) = q(y_1) + p(y) + B(y)|y_1|^{\gamma-1}\|y\|^\beta \qquad \forall y \in \mathbb{R}^d.$$

LEMMA 3. *Let $\mathcal{A} = (\gamma, \beta) \in \mathcal{P}_3$. Then, for any $\lambda > 0$ small enough,*

(26) $$\sup_{(q,p,B) \in \mathfrak{B}(\mathcal{A},\lambda)} |[\mathsf{K}_{(\mathcal{A},\lambda)} * Q](0) - Q(0)| \leq c\lambda,$$

(27) $$\int |\mathsf{K}_{(\mathcal{A},\lambda)}(y)|\|y\|^m\,du \leq c'\lambda^{m/(\gamma\beta)} \qquad \forall m \in \mathbb{R},$$

where the constant $c$ depends only on $c_3, c_4, L, d$ and $\mathcal{A}$, and $c'$ depends only on $m, d$ and $\mathcal{A}$. Furthermore,

$$\text{(28)} \qquad \|\mathsf{K}_{(\mathcal{A},\lambda)}\|_1 \leq c'' \quad \text{and} \quad \|\mathsf{K}_{(\mathcal{A},\lambda)}\|_2 \leq c^{(3)} \lambda^{-d/(2\beta)},$$

where the constants $c''$ and $c^{(3)}$ only depend on $\mathcal{A}$ and $d$.

The weight $\mathsf{K}_{(\mathcal{A},\lambda)}$ depends on $\mathcal{A} = (\gamma, \beta)$ in such a way that the constants in the bounds (26)–(28) diverge when $\mathcal{A}$ approaches the boundary $d(\gamma - 1) + 1 = \beta$ of the zone $\mathcal{P}_3$. So, Lemma 3 cannot be extended to $\mathcal{A} \in \mathcal{P}_2$.

We consider now another construction that provides the weight $\mathsf{K}_{(\mathcal{A},\lambda)}$ with the properties similar to those of Lemma 3 but satisfied for all $\mathcal{A} \in \mathcal{P}_2 \cup \mathcal{P}_3$ and, what is more, uniformly over this set. The price to pay for the uniformity is an extra $\log \log(1/\lambda)$ factor in the bound for the $L_1$-norm of $\mathsf{K}_{(\mathcal{A},\lambda)}$.

LEMMA 4. *Let $\mathcal{A} = (\gamma, \beta) \in \mathcal{P}_2 \cup \mathcal{P}_3$. Then, for any $\lambda > 0$ small enough,*

$$\text{(29)} \qquad \sup_{(q,p,B) \in \mathfrak{B}(\mathcal{A},\lambda)} |[\mathsf{K}_{(\mathcal{A},\lambda)} * Q](0) - Q(0)| \leq c_5 \lambda,$$

$$\text{(30)} \qquad \int |\mathsf{K}_{(\mathcal{A},\lambda)}(y)| \|y\|^m \, du \leq c_6 \lambda^{m/(\gamma\beta)} \qquad \forall m \in \mathbb{R},$$

where the constant $c_5$ depends only on $c_3, c_4, L$ and $d$, and $c_6 > 0$ depends only on $m$ and $d$ (both constants are explicit in the proof of the lemma). Furthermore,

$$\text{(31)} \quad \|\mathsf{K}_{(\mathcal{A},\lambda)}\|_1 \leq c_7 \ln \ln \lambda^{-1} \quad \text{and} \quad \|\mathsf{K}_{(\mathcal{A},\lambda)}\|_2 \leq c_8 \lambda^{-(\beta+d-1)/(2\gamma\beta)},$$

where the constants $c_7$ and $c_8$ only depend on $d$.

### 6.2. Weight systems.

*Weight system for zone $\mathcal{P}_4$ (no local structure).* The construction of $\mathsf{K}_{(\mathcal{A},\lambda)}$ is trivial when $\mathcal{A}$ is in the zone $\mathcal{P}_4$ of no local structure. In this case a basic boxcar kernel tuned to the smoothness of the composite function can be used. Observe that when $\mathcal{A} \in (0,1]^2$ the smoothness of the composite function equals to $\gamma\beta$, and when $\mathcal{A} = (\gamma, \beta)$ satisfies $1 < \beta \leq \gamma \leq 2$ the smoothness is $\beta$. So, we define the weight $\mathsf{K}_{(\mathcal{A},\lambda)}$ for the zone $\mathcal{P}_4$ as follows:

$$\mathsf{K}_{(\mathcal{A},\lambda)}(y) = \begin{cases} (2\lambda^{1/(\gamma\beta)})^{-d} \mathbb{I}_{[-\lambda^{1/(\gamma\beta)}, \lambda^{1/(\gamma\beta)}]^d}(y), & \text{if } \mathcal{A} = (\gamma, \beta) \in (0,1]^2, \\ (2\lambda^{1/\beta})^{-d} \mathbb{I}_{[-\lambda^{1/\beta}, \lambda^{1/\beta}]^d}(y), & \text{if } 1 < \beta < \gamma \leq 2. \end{cases}$$

Here $\mathbb{I}_A(\cdot)$ stands for the indicator function of a set $A$. The proof of Lemma 1 is straightforward.

*Weight system for zone* $\mathcal{P}_1$ *(local single-index model).* The zone of local single-index model is $\mathcal{P}_1 = \{\mathcal{A} = (\gamma, \beta) : \gamma \leq 1, 1 < \beta \leq 2\}$. For any $\mathcal{A} \in \mathcal{P}_1$ and $\lambda > 0$ consider the hyperrectangle

$$\Pi_\lambda(\mathcal{A}) = [-\lambda^{1/\gamma}, \lambda^{1/\gamma}] \times [-\lambda^{1/(\gamma\beta)}, \lambda^{1/(\gamma\beta)}]^{d-1}$$

and define the weight $\mathsf{K}_{(\mathcal{A},\lambda)}$ as follows:

(32) $$\mathsf{K}_{(\mathcal{A},\lambda)} = (2^d \lambda^{1/\gamma + (d-1)/(\gamma\beta)})^{-1} \mathbb{I}_{\Pi_\lambda(\mathcal{A})}(y), \qquad y \in \mathbb{R}^d.$$

The proof of Lemma 2 is evident.

*Weight system for zone* $\mathcal{P}_2 \cup \mathcal{P}_3$ *(RISD local model).* The zone of RISD local model is $\mathcal{P}_2 \cup \mathcal{P}_3 = \{\mathcal{A} = (\gamma, \beta) : 1 < \gamma \leq \beta \leq 2\}$. The definition of the weight in this case is more involved. Indeed, taking $\mathsf{K}_{(\mathcal{A},\lambda)}$ as a simple product of boxcar kernels (32) results for $\mathcal{A} \in \mathcal{P}_2 \cup \mathcal{P}_3$ in too large approximation error.

Our aim is to construct a weight $\mathsf{K}_{(\mathcal{A},\lambda)} : \mathbb{R}^d \to \mathbb{R}$ with the following properties:

– for some $c > 0$, it should vanish outside the set [cf. (8)]

$$\{y \in \mathbb{R}^d : |y_1| \leq c\lambda^{1/\beta}, \|y\| \leq c\lambda^{1/(\gamma\beta)}, |y_1|^{\gamma-1} \|y\|^\beta \leq c\lambda\}.$$

– for a function $q(y_1)$ of the first component $y_1$ of $y \in \mathbb{R}^d$, the "characteristic size" of $\mathsf{K}_{(\mathcal{A},\lambda)}$ should be $\lambda^{1/\gamma}$; for a function $Q(y_2, \ldots, y_d)$ of the remaining components $y_2, \ldots, y_d$ it should be $\lambda^{1/\beta}$. Namely, we want to ensure the relations

$$\int \mathsf{K}_{(\mathcal{A},\lambda)}(y) q(y_1) \, dy = (2\lambda^{1/\gamma})^{-1} \int_{-\lambda^{1/\gamma}}^{\lambda^{1/\gamma}} q(y_1) \, dy_1$$

and

$$\int \mathsf{K}_{(\mathcal{A},\lambda)}(y) Q(y_2, \ldots, y_d) \, dy$$
$$= (2\lambda^{1/\beta})^{-(d-1)} \int_{-\lambda^{1/b}}^{\lambda^{1/\beta}} \cdots \int_{-\lambda^{1/\beta}}^{\lambda^{1/\beta}} Q(y_2, \ldots, y_d) \, dy_2 \cdots dy_d.$$

These properties are crucial to guarantee that the bias of linear approximation is of the order $O(\lambda)$ (cf. Lemma 3). Note that the simple rectangular kernel (32) used for the local single-index model can attain such a bias, but only at the price of too large $L_2$-norm (which characterizes the variance). We now give an example showing how a weight with the required properties can be constructed in a particular case.

*The two-step weight.* Set

(33) $\quad u_1 = \lambda^{1/\gamma}, \quad u_2 = \lambda^{1/\beta}, \quad v_1 = \lambda^{(\beta-\gamma+1)/\beta^2}, \quad v_2 = \frac{1}{2}\lambda^{1/\beta},$

$$\Pi_{1,1} = [0, u_1] \times [v_2, v_1]^{d-1}, \qquad \mu_{1,1} = u_1(v_1 - v_2)^{d-1};$$
$$\Pi_{2,2} = [u_1, u_2] \times [0, v_2]^{d-1}, \qquad \mu_{2,2} = (u_2 - u_1)v_2^{d-1};$$
$$\Pi_{2,1} = [u_1, u_2] \times [v_2, v_1]^{d-1}, \qquad \mu_{2,1} = (u_2 - u_1)(v_1 - v_2)^{d-1}.$$

Next, we define, for $y \in \mathbb{R}_+^d$,

(34) $\quad \Lambda(y) = \mu_{1,1}^{-1}\mathbb{I}_{\Pi_{1,1}}(y) - \mu_{2,1}^{-1}\mathbb{I}_{\Pi_{2,1}}(y) + \mu_{2,2}^{-1}\mathbb{I}_{\Pi_{2,2}}(y).$

For $y = (y_1, \ldots, y_d) \in \mathbb{R}^d$ we write $|y| = (|y_1|, \ldots, |y_d|)$ and define the weight $\mathsf{K}_{(\mathcal{A},\lambda)}$ for $y \in \mathbb{R}^d$ by the relation

(35) $\quad \mathsf{K}_{(\mathcal{A},\lambda)}(y) = 2^{-d}\Lambda(|y|).$

We will call this weight the *two-step weight* (cf. Figure 4). Its key property is as follows. First, for any integrable function $q(y_1)$ of the first coordinate $y_1$ we have

$$\int \mathsf{K}_{(\mathcal{A},\lambda)}(y) q(y_1) \, dy = \frac{1}{2u_1} \int_{-u_1}^{u_1} q(y_1) \, dy_1,$$

since the integral of $q$ over $\Pi_{2,1}$ is exactly the same as that over $\Pi_{2,2}$. Further, for any integrable function $Q(y_2, \ldots, y_d)$ of $y_2, \ldots, y_d$,

$$\int \mathsf{K}_{(\mathcal{A},\lambda)}(y) Q(y_2, \ldots, y_d) \, dy$$
$$= (2v_2)^{-(d-1)} \int_{-v_2}^{v_2} \cdots \int_{-v_2}^{v_2} Q(y_2, \ldots, y_d) \, dy_2 \cdots dy_d,$$

since the integral of $Q$ over $\Pi_{2,1}$ is exactly the same as that over $\Pi_{1,1}$. In other words, the negative term $-\mu_{2,1}^{-1}\mathbb{I}_{\Pi_{2,1}}(y)$ in (34) allows us to compensate

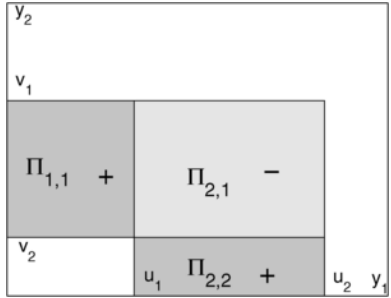

FIG. 4. *Pavement $\Pi_{i,j}$ for the two-step weight, $d = 2$. The weight vanishes in the white zones.*

26    A. B. JUDITSKY, O. V. LEPSKI AND A. B. TSYBAKOV

the excess of the bias introduced by the two other terms, so that the resulting bias remains of the order $O(\lambda)$ (cf. Lemma 3).

For the two-step weight (35) we have

$$\int \mathsf{K}_{(\mathcal{A},\lambda)}(y)\,dy = 1, \qquad \|\mathsf{K}_{(\mathcal{A},\lambda)}\|_1 = 3, \qquad \|\mathsf{K}_{(\mathcal{A},\lambda)}\|_2^2 = \mu_{1,1}^{-1} + \mu_{2,2}^{-1} + \mu_{2,1}^{-1}.$$

We now define

$$\rho = \frac{(d-1)(\gamma-1)}{\beta}$$

and consider the subset $\{\mathcal{A} = (\gamma,\beta) : \rho \geq (\beta-\gamma)/\gamma\}$ of $\mathcal{P}_3$. It is easy to see that for $\rho \geq (\beta-\gamma)/\gamma$ we have

$$\|\mathsf{K}_{(\mathcal{A},\lambda)}\|_2^2 = O(\lambda^{-d/\beta}).$$

Since $\gamma \leq \beta$ for $\mathcal{A} \in \mathcal{P}_3$, this result is better than part (ii) of Lemma 2 where $\mathsf{K}_{(\mathcal{A},\lambda)}$ is a rectangular kernel. But we need the condition $\rho \geq (\beta-\gamma)/\gamma$. It is clearly satisfied when $\rho \geq 1$ (recall that $\gamma > 1, \beta \leq 2$). For smaller values of $\rho$ we need to add extra "steps" in the construction, that is, to introduce piecewise constant weights with more and more pieces of the pavement, in order to get the bias compensation property as discussed above. For instance, if $\rho + \rho^2 \geq \frac{\beta-\gamma}{\gamma}$ [since $(\beta-\gamma)/\gamma < 1$, this is certainly the case when $\rho \geq \frac{\sqrt{5}-1}{2}$] we need a pavement of five sets $\Pi_{i,j}$ in order to obtain a piecewise constant weight with the required statistical properties, and so on. We come to the following construction of the weight.

*Generic construction.* Define a piecewise constant weight $\mathsf{K}_{(\mathcal{A},\lambda)}$ as follows. Fix an integer $r$ that we will further call *number of steps* (of weight construction). Let $(u_j)_{j=1,\ldots,r}$ and $(v_j)_{j=1,\ldots,r+1}$ be, respectively, a monotone increasing and a monotone decreasing sequence of positive numbers with $u_1 = \lambda^{1/\gamma}$, $v_r = \lambda^{1/\beta}/2$ and $v_{r+1} = 0$. We set

$$\Pi_{1,1} = [0, u_1] \times [v_2, v_1]^{d-1}, \qquad \mu_{1,1} = u_1(v_1 - v_2)^{d-1}.$$

For $i = 2, \ldots, r$ and $j = i-1, i$ we define

$$\Pi_{i,j} = [u_{i-1}, u_i] \times [v_{j+1}, v_j]^{d-1}, \qquad \mu_{i,j} = (u_i - u_{i-1})(v_j - v_{j+1})^{d-1}.$$

For $y \in \mathbb{R}_+^d$ consider

$$\Lambda_1(y) = \frac{1}{\mu_{1,1}} \mathbb{I}_{\Pi_{1,1}}(y);$$

$$\Lambda_i(y) = \frac{1}{\mu_{i,i}} \mathbb{I}_{\Pi_{i,i}}(y) - \frac{1}{\mu_{i,i-1}} \mathbb{I}_{\Pi_{i,i-1}}(y), \qquad i = 2, \ldots, r.$$

The weight $\mathsf{K}_{(\mathcal{A},\lambda)}$ is defined for $y = (y_1, \ldots, y_d) \in \mathbb{R}^d$ as follows:

$$\mathsf{K}_{(\mathcal{A},\lambda)}(y) = 2^{-d} \sum_{i=1}^{r} \Lambda_i(|y|), \tag{36}$$

where $|y| = (|y_1|, \ldots, |y_d|)$. Clearly,

$$\int \mathsf{K}_{(\mathcal{A},\lambda)}(y) \, dy = 1, \qquad \|\mathsf{K}_{(\mathcal{A},\lambda)}\|_1 = 2r - 1.$$

*Construction of the weight for* $\mathcal{A} \in \mathcal{P}_3 = \{\mathcal{A} : 1 < \gamma \leq \beta \leq 2, \beta < d(\gamma - 1) + 1\}$. If $\rho \geq \frac{\beta - \gamma}{\gamma}$ we define $\mathsf{K}_{(\mathcal{A},\lambda)}$ as a two-step weight, that is, we set $r = 2$ and take $(u_j)$ and $(v_j)$ as in (33).

If $\rho < \frac{\beta - \gamma}{\gamma}$ we use another definition. We introduce the sequence $(\alpha_k)_{k \geq 0}$ as follows:

$$\alpha_0 = \beta^{-1}, \qquad \alpha_{k+1} = \alpha_k \rho + \beta^{-1} = \beta^{-1} \sum_{i=0}^{k+1} \rho^i, \qquad k = 1, 2, \ldots. \tag{37}$$

The sequence $(\alpha_k)$ is monotone increasing and, since $\beta < d(\gamma - 1) + 1$, we have

$$\begin{aligned} &\lim_{k \to \infty} \alpha_k = \infty, \qquad \text{if } \rho \geq 1, \\ &\lim_{k \to \infty} \alpha_k = (\beta - (\gamma - 1)(d - 1))^{-1} > \frac{1}{\gamma}, \qquad \text{if } \rho < 1. \end{aligned} \tag{38}$$

Thus we can define an integer $r \geq 2$ such that

$$\alpha_{r-1} \geq \frac{1}{\gamma} > \alpha_{r-2}. \tag{39}$$

Note that $r$ depends only on $\mathcal{A} = (\gamma, \beta)$ and $d$. Now we set

$$\begin{aligned} u_1 &= \lambda^{1/\gamma}, \qquad u_i = \lambda^{\alpha_{r-i}}, \qquad i = 2, \ldots, r; \\ v_i &= \lambda^{1/\beta} u_{i+1}^{-(\gamma-1)/\beta}, \qquad i = 1, \ldots, r - 1. \end{aligned} \tag{40}$$

Recall that $v_r = \frac{1}{2}\lambda^{1/\beta}$ and $v_{r+1} = 0$. If $\rho < \frac{\beta - \gamma}{\gamma}$ define the weight $\mathsf{K}_{(\mathcal{A},\lambda)}$ by (36), with the sequences $(u_j)$ and $(v_j)$ as in (40).

Note that for $\rho \geq \frac{\beta - \gamma}{\gamma}$ the weight $\mathsf{K}_{(\mathcal{A},\lambda)}$ is just the two-step weight. The corresponding pavement $\{\Pi_{i,j}\}$ only contains three sets (cf. Figure 4).

*Construction of the weight for* $\mathcal{A} \in \mathcal{P}_2$. We consider now another choice of the sequences $(u_i)$ and $(v_i)$, which provides the weight $\mathsf{K}_{(\mathcal{A},\lambda)}$ with the properties similar to those of Lemma 3 but satisfied for all $\mathcal{A} \in \mathcal{P}_2 \cup \mathcal{P}_3$ and, what is more, uniformly over this set. The price to pay for the uniformity is an extra $\log\log(1/\lambda)$ factor in the bound for the $L_1$-norm of $\mathsf{K}_{(\mathcal{A},\lambda)}$.

If $(\beta-\gamma)/\gamma \leq (1+\rho)\rho$ we define the weight as in Lemma 3. If $(\beta-\gamma)/\gamma > (1+\rho)\rho$ we use another definition of sequences $(u_i)$ and $(v_i)$. For any $0 < \lambda < 1$ we define

$$
(41) \qquad V(\lambda) = \ln\left\{\frac{(\gamma-1)(\beta-\gamma)}{\gamma\beta^2}\ln(1/\lambda)\right\}.
$$

If $V(\lambda) \leq 0$ we define $\mathsf{K}_{(\mathcal{A},\lambda)}$ as a two-step weight, that is, we set $r=2$ and take $(u_j)$ and $(v_j)$ as in (33). If $V(\lambda) > 0$ we define $r = r(\lambda) > 1$ by

$$
r = \min\left\{s \in \mathbb{N} : s > 1, \frac{V(\lambda)}{s-1} < \frac{1}{2}\ln\left(\frac{\sqrt{5}+1}{2}\right)\right\}.
$$

Next, set $\alpha = \frac{V(\lambda)}{r-1}$, $\nu = (\frac{\sqrt{5}+1}{2})^{1/2}$ and define the sequences $(u_i)$ and $(v_i)$ as follows

$$
\begin{aligned}
u_i &= \lambda^{1/\gamma}\exp\left\{\frac{\beta}{\gamma-1}\exp(\alpha(i-1))\right\}, \qquad i=1,\ldots,r,\\
v_i &= \lambda^{1/(\gamma\beta)}\exp\{-\nu\exp(\alpha i)\}, \qquad i=1,\ldots,r-1, \qquad v_r = \tfrac{1}{2}\lambda^{1/\beta}.
\end{aligned}
$$
(42)

Note that $u_r = \lambda^{1/\beta}$.

Some remarks are in order here.

1. The number of steps $r$ in the construction of the weight is typically small. In particular, $r=2$ if $\rho \geq \frac{\beta-\gamma}{\gamma}$, and $r=3$ if $(1+\rho)\rho \geq \frac{\beta-\gamma}{\gamma} > \rho$ [cf. (39)]. Moreover, for $1 < \gamma \leq \beta \leq 2$ we have

$$
\frac{(\gamma-1)(\beta-1)}{\gamma\beta^2} \leq \frac{(\beta-1)^2}{\beta^3} \leq \frac{1}{8}.
$$

Hence, $V(\lambda) \leq \ln(\frac{\sqrt{5}+1}{2})$ for all $\lambda > 3 \cdot 10^{-6}$, which means that for $(1+\rho)\rho < \frac{\beta-\gamma}{\gamma}$ no more than 3 steps of the construction are needed if $\lambda > 3 \cdot 10^{-6}$. In other words, unless we are not "extremely far" in the asymptotics, the number of steps $r$ does not exceed 3 and thus the $L_1$-norm of the resulting weight $\mathsf{K}_{(\mathcal{A},\lambda)}$ is bounded by 5.

2. In the asymptotics when $\lambda \to 0$ the number of steps $r = r(\lambda)$ in the construction and thus the $L_1$-norm of the weight $\mathsf{K}_{(\mathcal{A},\lambda)}$ is at most $O(\ln\ln\lambda^{-1})$. As discussed in the previous remark, this behavior starts "extremely far" in the asymptotics, so it has essentially a theoretical interest. In the theory, it results in an extra $\ln\ln\varepsilon^{-1}$ factor in the upper bound for the estimation procedure, as compared to the lower bound in (5). It can be shown that for $\mathcal{A} \in \mathcal{P}_2$ a weight with the required approximation properties cannot have the $L_1$-norm growing slower than $\ln\ln\lambda^{-1}$, as $\lambda \to 0$. On the other hand, as we have seen in Lemma 3, for $\mathcal{A} \in \mathcal{P}_3$ solely, there is a choice of sequences $(u_j)$ and $(v_j)$ such that the $L_1$-norm of the weight is bounded by a constant independent of $\lambda$. This constant, however, depends on $\mathcal{A} = (\gamma, \beta)$ and explodes as $\mathcal{A}$ approaches the boundary of $\mathcal{P}_3$.

## 7. Proofs.

7.1. *Proof of Theorem 1.* For any $\beta > 0, \gamma > 0$ and any $0 < \varepsilon < 1$ define the integers

$$q_1 = \lceil (\varepsilon \sqrt{\ln(1/\varepsilon)})^{-2/(2\gamma\beta + \beta + (d-1))} \rceil.$$

Consider the regular grid $\Gamma_{q_1}$ on $[0,1]^{d-1}$ defined by

$$\Gamma_{q_1} \triangleq \left\{ \left( \frac{2k_1 + 1}{2q_1}, \ldots, \frac{2k_{d-1} + 1}{2q_1} \right) : k_i \in \{0, \ldots, q_1 - 1\}, i = 1, \ldots, d-1 \right\}.$$

Denote by $x_1, \ldots, x_m$, where $m = \mathrm{card}(\Gamma_{q_1}) = q_1^{d-1}$, the elements of $\Gamma_{q_1}$ numbered in an arbitrary order.

Let $f_0 : \mathbb{R} \to \mathbb{R}_+$ be an infinitely differentiable function such that $f_0(0) = 1$, $f_0(u) = f_0(-u)$ for all $u \in \mathbb{R}$, $f_0(u) = 0$ for $u \notin [-1/2, 1/2]$, and $f_0$ is strictly monotone decreasing on $[0, 1/2]$. Examples of such functions can be readily constructed; compare [27], page 78. Set

$$\varphi_0(t_2, \ldots, t_d) = \frac{1}{2} \prod_{j=2}^d f_0(t_j) \qquad \forall (t_2, \ldots, t_d) \in \mathbb{R}^{d-1}$$

and

$$f(u) = L_0 h^\gamma f_0 \left( \frac{u}{h} \right) \qquad \forall u \in \mathbb{R},$$

where $h = h_1^\beta$, $h_1 = 1/q_1$ and $0 < L_0 < 1$ is a constant to be chosen small enough. Consider the following collection of infinitely differentiable functions of $t = (t_1, \ldots, t_d) \in \mathbb{R}^d$:

$$g_k(t) = f(G_k(t)) = L_0 h^\gamma f_0 \left( \frac{G_k(t)}{h} \right), \qquad k = 0, 1, \ldots, m,$$

where

$$G_0(t) = L_0 \sin t_1,$$

$$G_k(t) = L_0 \sin t_1 + L_0 h_1^\beta \varphi_0 \left( \frac{t_2 - x_{k,2}}{h_1}, \ldots, \frac{t_d - x_{k,d}}{h_1} \right), \qquad k = 1, \ldots, m$$

and $x_{k,j}$ stands for the $j$th component of $x_k$. We note that, in view of the above definitions, the sets where the functions $g_l$ and $g_k$ differ from $g_0$ are disjoint for $l \neq k$, $k \neq 0$, $l \neq 0$.

It is easy to see that if $L_0$ is small enough, $g_k \in \mathbb{H}(\mathcal{A}, \mathcal{L})$, $k = 0, \ldots, m$. In what follows, we assume that $L_0$ is chosen in this way. To prove Theorem 1, we follow the scheme of lower bounds based on reduction to the problem of testing $m+1$ hypotheses (cf., e.g., [27]). We choose the hypotheses to be determined by $g_0, \ldots, g_m$ and we apply Theorem 2.5 of [27], where we consider the sup-norm distance $d(g_l, g_k) = \|g_l - g_k\|_\infty = \sup_{t \in [-1,1]^d} |g_l(t) - g_k(t)|$, $l, k = 0, 1, \ldots, m$. Since the functions $g_l$ and $g_k$ differ from $g_0$ on disjoint sets, for any $l \neq k, l, k = 1, \ldots, m$, we have

$$d(g_l, g_k) = d(g_0, g_k) \geq L_0 h^\gamma |f_0(0) - f_0(L_0 h_1^\beta \varphi_0(0)/h)|$$
$$= L_0 h^\gamma |f_0(0) - f_0(L_0(1 + o_\varepsilon(1))/2)|,$$

where $o_\varepsilon(1) \to 0$, as $\varepsilon \to 0$. Since $L_0 > 0$ and $f_0$ is strictly decreasing on $[0, \infty)$ there exists a constant $L^* > 0$ such that, for $\varepsilon$ small enough,

$$(43) \qquad d(g_l, g_k) \geq L^* h^\gamma \asymp (\varepsilon\sqrt{\ln(1/\varepsilon)})^{(2\gamma)/(2\gamma+1+(d-1)/\beta)},$$

$$l \neq k, l, k = 0, \ldots, m.$$

Thus, assumption (i) of Theorem 2.5 in [27] is satisfied with $s = L^* h^\gamma/2$. It remains to check assumption (ii) of that theorem. The probability measures $\mathbb{P}_{g_k}$ are Gaussian, and the Kullback–Leibler divergence between $\mathbb{P}_{g_k}$ and $\mathbb{P}_{g_0}$ has the form

$$\mathbf{K}(\mathbb{P}_{g_k}, \mathbb{P}_{g_0}) = \varepsilon^{-2} \int_\mathcal{D} (g_0(t) - g_k(t))^2 \, dt$$
$$= \varepsilon^{-2} L_0^2 h^{2\gamma} \int_\mathcal{D} \left| f_0\left(\frac{L_0 \sin t_1}{h}\right) - f_0\left(\frac{L_0 \sin t_1}{h} + w(t_2, \ldots, t_d)\right) \right|^2 dt,$$

where we write for brevity

$$w(t_2, \ldots, t_d) \triangleq L_0 \varphi_0\left(\frac{t_2 - x_{k,2}}{h_1}, \ldots, \frac{t_d - x_{k,d}}{h_1}\right).$$

Since, for any $a, w \in \mathbb{R}$,

$$\left| f_0\left(\frac{a}{h}\right) - f_0\left(\frac{a}{h} + w\right) \right|^2 = w^2 \left| \int_0^1 f_0'\left(\frac{a}{h} + uw\right) du \right|^2$$
$$\leq w^2 \int_0^1 \left| f_0'\left(\frac{a}{h} + uw\right) \right|^2 du,$$

we find

$$\mathbf{K}(\mathbb{P}_{g_k}, \mathbb{P}_{g_0})$$
$$\leq \varepsilon^{-2} L_0^2 h^{2\gamma} \int w^2(t_2, \ldots, t_d) \, dt_2 \cdots dt_d$$
$$\times \int_0^1 \left[ \int_{|t_1| \leq |\mathcal{D}|} \left| f_0'\left(\frac{L_0 \sin t_1}{h} + uw(t_2, \ldots, t_d)\right) \right|^2 dt_1 \right] du,$$

where $|\mathcal{D}|$ is the Euclidean diameter of $\mathcal{D}$. Since $f_0$ is supported on $[-1/2, 1/2]$ and $|w(t_2, \ldots, t_d)| < 1/2$, the values $f_0'((L_0 \sin t_1)/h + u w(t_2, \ldots, t_d))$ under the last integral can be nonzero only if $L_0 |\sin t_1| \leq h$. The Lebesgue measure of the set $\{t_1 : |t_1| \leq |\mathcal{D}|, L_0 |\sin t_1| \leq h\}$ is $O(h)$, as $h \to 0$. Hence, the double integral in the last display is bounded by $c_* h$ for all $h$ small enough, where $c_* > 0$ is an absolute constant. This yields

$$\mathbf{K}(\mathbb{P}_{g_k}, \mathbb{P}_{g_0}) \leq c_* L_0^4 \varepsilon^{-2} h^{2\gamma+1} h_1^{d-1} \int_{\mathbb{R}^{d-1}} \varphi_0^2(v)\, dv$$

$$\leq c_{**} L_0^4 \ln(1/\varepsilon),$$

where $c_{**} > 0$ is an absolute constant. Next, $m = q_1^{d-1}$, so that $\ln m \asymp \ln(1/\varepsilon)$. This and the previous inequality imply that if $L_0$ is chosen small enough, we have

(44) $$\mathbf{K}(\mathbb{P}_{g_k}, \mathbb{P}_{g_0}) \leq (1/16) \ln m.$$

Using (43), (44) and applying Theorem 2.5 in [27] we get the lower bound

(45) $$\liminf_{\varepsilon \to 0} \inf_{\tilde{g}_\varepsilon} \sup_{g \in \mathbb{H}(\mathcal{A}, \mathcal{L})} \mathbb{E}_g[((\varepsilon\sqrt{\ln(1/\varepsilon)})^{-(2\gamma)/(2\gamma+1+(d-1)/\beta)} \times \|\tilde{g}_\varepsilon - g\|_\infty)^p] > 0,$$

which is valid for all $\beta > 0, \gamma > 0$ and all $p > 0$.

We now show that for the trivial cases discussed in Section 2 we can obtain better lower bounds. Consider first the case where $0 < \beta, \gamma \leq 1$. Then we use the same technique as above, but we set now $q_1 = \lceil (\varepsilon\sqrt{\ln(1/\varepsilon)})^{-2/(2\gamma\beta+d)} \rceil$. We then introduce a regular grid $\Gamma_{q_1}^*$ on $[0,1]^d$ defined by

$$\Gamma_{q_1}^* \triangleq \left\{ \left(\frac{2k_1+1}{2q_1}, \ldots, \frac{2k_d+1}{2q_1}\right) : k_i \in \{0, \ldots, q_1-1\}, i = 1, \ldots, d \right\}$$

and denote by $x_1, \ldots, x_m$, where $m = \text{card}(\Gamma_{q_1}^*) = q_1^d$, the elements of $\Gamma_{q_1}^*$ numbered in an arbitrary order. We set now

$$\varphi_0(t) \triangleq \prod_{j=1}^d u(t_j) \qquad \forall t \in \mathbb{R}^d$$

and we choose the functions $g_k$ in the following way:

$$g_0(t) \equiv 0,$$

$$g_k(t) = \left| L_0 h^\beta \varphi_0\left(\frac{t - x_k}{h}\right) \right|^\gamma, \qquad t \in \mathbb{R}^d, \qquad k = 1, \ldots, m,$$

where $h = 1/q_1$. Note that for sufficiently small $h$ we can write these functions as compositions $g_k = f \circ G_k$, where $f(u) = L_0' |u|^\gamma f_0(u)$, $G_0(t) \equiv 0$, $G_k(t) = L_0' h^\beta \varphi_0((t-x_k)/h)$ and $L_0' = L_0^{\gamma/(\gamma+1)}$ with a slightly different definition of $f_0$ than above. Namely, we choose $f_0$ to be infinitely differentiable, supported on $[-1/2, 1/2]$ and such that $f_0(u) = 1$ for $u \in [-1/4, 1/4]$. It is easy to see that if $L_0$ is small enough, $g_k \in \mathbb{H}(\mathcal{A}, \mathcal{L})$, $k = 0, \ldots, m$. With this choice of $g_k$ we get

$$d(g_l, g_k) \geq L_0^\gamma h^{\gamma\beta} \varphi_0^\gamma(0) \asymp (\varepsilon\sqrt{\ln(1/\varepsilon)})^{(2\gamma\beta)/(2\gamma\beta+d)}, \tag{46}$$

$$l \neq k, \; l, \; k = 0, \ldots, m.$$

Next,

$$\mathbf{K}(\mathbb{P}_{g_k}, \mathbb{P}_{g_0}) = \varepsilon^{-2} \int_\mathcal{D} (g_0(t) - g_k(t))^2 \, dt$$

$$\leq L_0^{2\gamma} \varepsilon^{-2} h^{2\gamma\beta+d} \int_{\mathbb{R}^d} \varphi_0^{2\gamma}(v) \, dv \tag{47}$$

$$= O(\ln(1/\varepsilon)) \qquad \text{as } \varepsilon \to 0.$$

Using (46), (47) and Theorem 2.5 in [27], the proof is completed as in the previous case, so that we get the lower bound

$$\liminf_{\varepsilon \to 0} \inf_{\tilde{g}_\varepsilon} \sup_{g \in \mathbb{H}(\mathcal{A}, \mathcal{L})} \mathbb{E}_g[((\varepsilon\sqrt{\ln(1/\varepsilon)})^{-(2\gamma\beta)/(2\gamma\beta+d)} \|\tilde{g}_\varepsilon - g\|_\infty)^p] > 0, \tag{48}$$

which is valid for all $0 < \beta, \gamma \leq 1$ and all $p > 0$.

Finally, the second trivial case where (45) can be improved corresponds to $\gamma \geq \beta \vee 1$. As observed in Section 2, in this case we have the inclusion $\mathbb{H}_d(\beta, L_4) \subset \mathbb{H}(\mathcal{A}, \mathcal{L})$ with some constant $L_4 > 0$, and we can use the standard lower bound for $\mathbb{H}_d(\beta, L_4)$ (cf. [2, 3, 6, 23]):

$$\liminf_{\varepsilon \to 0} \inf_{\tilde{g}_\varepsilon} \sup_{g \in \mathbb{H}(\mathcal{A}, \mathcal{L})} \mathbb{E}_g[((\varepsilon\sqrt{\ln(1/\varepsilon)})^{-(2\beta)/(2\beta+d)} \|\tilde{g}_\varepsilon - g\|_\infty)^p] > 0. \tag{49}$$

Combining the bounds (45), (48) and (49) we obtain the result of Theorem 1.

7.2. *Proof of Theorem 2.* We need the following technical result.

LEMMA 5. *Let $\zeta = (\zeta_1, \ldots, \zeta_\mathcal{M})$ be a Gaussian random vector defined on a probability space $(\Omega, \mathcal{F}, \mathbb{P})$ and such that $\mathbb{E}\zeta_m = 0, \mathbb{E}\zeta_m^2 = \sigma_m^2, m = 1, \ldots, \mathcal{M}$. Let $\mathfrak{m}$ be a random variable with the values in $(1, \ldots, \mathcal{M})$ defined on the same probability space. Then for all $A > 1$ and all $s > 0$ we have*

$$\mathbb{E}(|\zeta_\mathfrak{m}|^s) \leq (\sqrt{2A\ln(\mathcal{M})})^s \left\{ \mathbb{E}(\sigma_\mathfrak{m}^s) + c_{12}(A,s)\mathcal{M}^{1-A} \max_{m=1,\ldots,\mathcal{M}} \sigma_m^s \right\},$$

*where $c_{12}(A, s) > 0$ is a constant depending only on $A$ and $s$.*

Proof is standard (see, e.g., [14]).

To prove Theorem 2 we proceed in steps.

*1°. Reduction to the discrete norm.* Fix $\mathcal{A} = (\gamma, \beta) \in (0, 2]^2$, and suppose that $g \in \mathbb{H}(\mathcal{A}, \mathcal{L})$. Let, for brevity, $\bar{g}_\varepsilon^* = g_{\mathcal{A}, \varepsilon}^*$. In view of the construction of the global estimator [cf. (19)] we get, for all $g \in \mathbb{H}(\mathcal{A}, \mathcal{L})$,

$$\|\bar{g}_\varepsilon^* - g\|_\infty \leq \sup_{z \in \mathbb{Z}^d} \max_{x \in \Pi_\varepsilon(z) \cap [-1, 1]^d} |\bar{g}_\varepsilon^*(x) - g(x)|$$

(50)
$$\leq |\bar{g}_\varepsilon^* - g|_\infty + C\varepsilon^{2\gamma(\beta \wedge 1)},$$

where

$$|\bar{g}_\varepsilon^* - g|_\infty \triangleq \max_{z \in \mathcal{Z}_\varepsilon} |\bar{g}_\varepsilon^*(z) - g(z)| \quad \text{with } \mathcal{Z}_\varepsilon = (\varepsilon^2 \mathbb{Z})^d \cap [-1, 1]^d.$$

Here and in what follows we will use the same notation $C$ for possibly different positive constants depending only on $\mathcal{A}, \mathcal{L}$ and $d$. Since $\varepsilon^{2\gamma(\beta \wedge 1)} = o(\phi_\varepsilon(\gamma, \beta)), \varepsilon \to 0$, for all $(\gamma, \beta) \in \mathbb{R}_+^2$, it is sufficient to prove Theorem 2 with the loss given by the maximum norm $|\cdot|_\infty$ on the finite set $\mathcal{Z}_\varepsilon$. *Thus, without loss of generality, in what follows we will replace $\|\cdot\|_\infty$ by $|\cdot|_\infty$.*

*2°. Control of large deviations.* To any $z \in \mathcal{Z}_\varepsilon$ we assign a vector $\theta^z \in \mathbb{S}_\varepsilon$ such that $\|\theta^z - \theta_0^z\| \leq \varepsilon^2$ where $\theta_0^z$ is defined in (20). Next, we set $\mathcal{J}_0^z \triangleq (\mathcal{A}, \theta^z, \lambda_\varepsilon(\mathcal{A}))$. Introduce the random event

$$\mathcal{F} = \{\exists z \in \mathcal{Z}_\varepsilon : \mathcal{J}_0^z \notin \hat{\mathfrak{T}}_z\},$$

where $\hat{\mathfrak{T}}_z$ is the set of acceptable triplets $\mathcal{J}$ defined in Section 5. We now show that for all $\varepsilon > 0$ small enough

(51) $$\sup_{g \in \mathbb{H}(\mathcal{A}, \mathcal{L})} \mathbb{P}_g(\mathcal{F}) \leq c_{12} \varepsilon^{2p},$$

where the constant $c_{12}$ depends only on $d$. Indeed, in view of the definition of the random set $\hat{\mathfrak{T}}_z$,

$$\mathcal{F} \subseteq \bigcup_{z \in \mathcal{Z}_\varepsilon} \bigcup_{\mathcal{J}' \in \mathfrak{J}_{\text{grid}}} \{|\Delta_{\mathcal{J}'} \hat{g}_{\mathcal{J}' * \mathcal{J}_0^z}(z)| > \mathbf{TH}_\varepsilon(\mathcal{J}', \mathcal{J}_0^z)\}$$

and therefore

(52) $$\mathbb{P}_g(\mathcal{F}) \leq \sum_{z \in \mathcal{Z}_\varepsilon} \sum_{\mathcal{J}' \in \mathfrak{J}_{\text{grid}}} \mathbb{P}_g\{|\Delta_{\mathcal{J}'} \hat{g}_{\mathcal{J}' * \mathcal{J}_0^z}(z)| > \mathbf{TH}_\varepsilon(\mathcal{J}', \mathcal{J}_0^z)\}.$$

Note that

$$\mathbb{E}_g \Delta_{\mathcal{J}'} \hat{g}_{\mathcal{J}' * \mathcal{J}_0^z}(z) = [\Delta_{\mathcal{J}'} K_{\mathcal{J}' * \mathcal{J}_0^z} * g](z).$$

Applying Proposition 2 with $\mathcal{J}_0^z = (\mathcal{A}, \theta^z, \lambda_\varepsilon(\mathcal{A}))$ and $\lambda = \lambda_0 = \lambda_\varepsilon(\mathcal{A})$ we obtain,

(53) $$\sup_{g \in \mathbb{H}(\mathcal{A}, \mathcal{L})} |\mathbb{E}_g \Delta_{\mathcal{J}'} \hat{g}_{\mathcal{J}' * \mathcal{J}_0^z}(z)|$$
$$\leq c_{11} \{\lambda_\varepsilon(\mathcal{A})(\|K_{\mathcal{J}'}\|_1 + \|K_{\mathcal{J}_0^z}\|_1) + \|K_{\mathcal{J}'}\|_1 \|K_{\mathcal{J}_0^z}\|_1 \varepsilon^2\}.$$

Now, due to the construction of the weight $\mathsf{K}_{(\mathcal{A},\lambda)}$ and the fact that $\|K_{\mathcal{J}}\|_1 = \|\mathsf{K}_{(\mathcal{A},\lambda_\varepsilon(\mathcal{A}))}\|_1$ for all $\mathcal{J} \in \mathfrak{J}_{\mathrm{grid}}$, there exists a constant $c_{13}$ depending only on $\mathcal{A}$ and $d$ such that $K_{\mathcal{A}}^* \triangleq \max_{\mathcal{J} \in \mathfrak{J}_{\mathrm{grid}}} \|K_{\mathcal{J}}\|_1$ satisfies

$$K_{\mathcal{A}}^* \leq c_{13}, \qquad \text{if } \mathcal{A} \in (0,2]^2 \setminus \mathcal{P}_2,$$
$$K_{\mathcal{A}}^* \leq c_{13} \ln\ln(1/\varepsilon), \qquad \text{if } \mathcal{A} \in \mathcal{P}_2.$$

Since also $\|K_{\mathcal{J}}\|_1 \geq 1$ and $\lambda_\varepsilon(\mathcal{A})/(\varepsilon \ln\ln(1/\varepsilon)) \to \infty$, as $\varepsilon \to 0$, we have, for $\varepsilon > 0$ small enough,

$$\sup_{g \in \mathbb{H}(\mathcal{A},\mathcal{L})} |\mathbb{E}_g \Delta_{\mathcal{J}'} \hat{g}_{\mathcal{J}' * \mathcal{J}_0^z}(z)|$$
(54)
$$\leq 2 c_{11} \lambda_\varepsilon(\mathcal{A})(\|K_{\mathcal{J}'}\|_1 + \|K_{\mathcal{J}_0^z}\|_1)$$
$$= 2\varepsilon \sqrt{\ln(1/\varepsilon)} \|\mathsf{K}_{(\mathcal{A},\lambda_\varepsilon(\mathcal{A}))}\|_2 (\|K_{\mathcal{J}'}\|_1 + \|K_{\mathcal{J}_0^z}\|_1),$$

where we used that $\lambda_\varepsilon(\mathcal{A})$ is a solution of (17). Note also that in $\mathbb{P}_g$-probability

(55)    $\Delta_{\mathcal{J}'} \hat{g}_{\mathcal{J}' * \mathcal{J}_0^z}(z) - \mathbb{E}_g \Delta_{\mathcal{J}'} \hat{g}_{\mathcal{J}' * \mathcal{J}_0^z}(z) \sim \mathcal{N}(0, \varepsilon^2 \|\Delta_{\mathcal{J}'} K_{\mathcal{J}' * \mathcal{J}_0^z}\|_2^2).$

Using (22), (52)–(55) and the definition of the threshold $\mathbf{TH}_\varepsilon(\cdot,\cdot)$ we obtain that, for $\varepsilon > 0$ small enough,

$$\mathbb{P}_g(\mathcal{F}) \leq \mathrm{card}(\mathcal{Z}_\varepsilon) \mathrm{card}(\mathbb{S}_\varepsilon) \mathbb{P}\{|\xi| > \sqrt{(4p+8d)\ln(1/\varepsilon)}\}$$
$$\leq \mathrm{card}(\mathcal{Z}_\varepsilon) \mathrm{card}(\mathbb{S}_\varepsilon) \varepsilon^{2p+4d},$$

where $\xi \sim \mathcal{N}(0,1)$. This proves (51) since $\mathrm{card}(\mathcal{Z}_\varepsilon) \leq (2\varepsilon^{-2}+1)^d$ and $\mathrm{card}(\mathbb{S}_\varepsilon) \leq (\sqrt{d}/\varepsilon)^{d-1}$.

$3°$. *Two intermediate bounds on the risks.* Using that $|\bar{g}_\varepsilon^*| \leq \ln\ln(1/\varepsilon)$ and $g \in \mathbb{H}(\mathcal{A},\mathcal{L})$ is uniformly bounded we deduce from (51) that, for all $\mathcal{A} = (\gamma,\beta) \in (0,2]^2$,

(56)    $\limsup_{\varepsilon \to 0} \sup_{g \in \mathbb{H}(\mathcal{A},\mathcal{L})} \mathbb{E}_g(\phi_\varepsilon^{-p}(\gamma,\beta)|\bar{g}_\varepsilon^* - g|_\infty^p \mathbb{I}\{\mathcal{F}\}) = 0.$

We now control the bias of $\hat{g}_{\mathcal{J}_0^z}$ via Proposition 1, its stochastic error via the bounds on $\|\mathsf{K}_{(\mathcal{A},\lambda_\varepsilon(\mathcal{A}))}\|_2$ in Lemmas 2–4 and apply (17) to get that, for all $\mathcal{A} = (\gamma,\beta) \in (0,2]^2$,

(57)    $\limsup_{\varepsilon \to 0} \sup_{g \in \mathbb{H}(\mathcal{A},\mathcal{L})} \mathbb{E}_g(\phi_\varepsilon^{-p}(\gamma,\beta)|\hat{g}_{\mathcal{J}_0^z} - g|_\infty^p) < \infty.$

$4°$. *Final argument.* Note that on the event $\mathcal{F}^c$ the set $\hat{\mathfrak{T}}_z$ of acceptable triplets $\mathcal{J}$ is nonempty for every $z \in \mathcal{Z}_\varepsilon$, so that $\hat{\mathcal{J}}_z$ exists. Thus, on $\mathcal{F}^c$ we can write, for all $z \in \mathcal{Z}_\varepsilon$,

$$(58) \quad |\hat{g}_{\hat{\mathcal{J}}_z}(z) - g(z)| \leq |\Delta_{\hat{\mathcal{J}}_z} \hat{g}_{\hat{\mathcal{J}}_z * \mathcal{J}_0^z}(z)| + |\Delta_{\mathcal{J}_0^z} \hat{g}_{\mathcal{J}_0^z * \hat{\mathcal{J}}_z}(z)| + |\hat{g}_{\mathcal{J}_0^z}(z) - g(z)|.$$

Further, on $\mathcal{F}^c$ the triplet $\mathcal{J}_0^z$ is acceptable for all $z \in \mathcal{Z}_\varepsilon$. This and the acceptability (by definition) of $\hat{\mathcal{J}}_z$ imply that on $\mathcal{F}^c$, for all $z \in \mathcal{Z}_\varepsilon$,

$$(59) \quad \begin{aligned} |\Delta_{\mathcal{J}_0^z} \hat{g}_{\mathcal{J}_0^z * \hat{\mathcal{J}}_z}(z)| &\leq \mathbf{TH}_\varepsilon(\mathcal{J}_0^z, \hat{\mathcal{J}}_z), \\ |\Delta_{\hat{\mathcal{J}}_z} \hat{g}_{\hat{\mathcal{J}}_z * \mathcal{J}_0^z}(z)| &\leq \mathbf{TH}_\varepsilon(\hat{\mathcal{J}}_z, \mathcal{J}_0^z). \end{aligned}$$

This, the definition of the threshold $\mathbf{TH}_\varepsilon$ and the fact that $\|K_\mathcal{J}\|_2 = \|\mathsf{K}_{(\mathcal{A},\lambda_\varepsilon(\mathcal{A}))}\|_2$ for all $\mathcal{J} \in \mathfrak{J}_{\text{grid}}$ yield that on $\mathcal{F}^c$, for all $z \in \mathcal{Z}_\varepsilon$,

$$(60) \quad \begin{aligned} |\hat{g}_{\hat{\mathcal{J}}_z}(z) - g(z)| &\leq 4C(p,d) K_\mathcal{A}^* \|\mathsf{K}_{(\mathcal{A},\lambda_\varepsilon(\mathcal{A}))}\|_2 \varepsilon \sqrt{\ln(1/\varepsilon)} + |\hat{g}_{\mathcal{J}_0^z}(z) - g(z)| \\ &= 4C(p,d) c_{11}^{-1} K_\mathcal{A}^* \lambda_\varepsilon(\mathcal{A}) + |\hat{g}_{\mathcal{J}_0^z}(z) - g(z)|. \end{aligned}$$

We combine (57) and (60) to get, with some constants $c_{14} - c_{16}$ independent of $\varepsilon$,

$$(61) \quad \begin{aligned} \sup_{g \in \mathbb{H}(\mathcal{A},\mathcal{L})} \mathbb{E}_g(|\bar{g}_\varepsilon^* - g|_\infty^p \mathbb{I}\{\mathcal{F}^c\}) &\leq c_{14}(K_\mathcal{A}^* \lambda_\varepsilon(\mathcal{A}))^p + c_{15} \phi_\varepsilon^p(\gamma,\beta) \\ &\leq c_{16}(K_\mathcal{A}^* \phi_\varepsilon(\gamma,\beta))^p. \end{aligned}$$

Theorem 2 follows now from (56) and (61).

## APPENDIX: PROOFS OF AUXILIARY RESULTS

### A.1. Proof of Proposition 2.

$1°$. PRELIMINARY REMARKS. For any $\mathcal{J} \in \mathfrak{J}$ and any $x \in [-1,1]^d$ we may write

$$\begin{aligned} [\Delta_\mathcal{J} K_{\mathcal{J} * \mathcal{J}_0^x} * g](x) &= [K_{\mathcal{J} * \mathcal{J}_0^x} * g](x) - [K_\mathcal{J} * g](x) \\ &= \int \left( \int K_\mathcal{J}(y-x) K_{\mathcal{J}_0^x}(t-y) \, dy \right) g(t) \, dt - [K_\mathcal{J} * g](x) \\ &= \int K_\mathcal{J}(y-x) \left( \int K_{\mathcal{J}_0^x}(t-y) g(t) \, dt \right) dy - [K_\mathcal{J} * g](x) \\ &= \int K_\mathcal{J}(y-x) g(y) \, dy - [K_\mathcal{J} * g](x) \end{aligned}$$

$$
\begin{aligned}
(62) \quad &+ \int K_{\mathcal{J}}(y-x)\left(\int K_{\mathcal{J}_0^x}(t-y)[g(t)-g(y)]\,dt\right)dy \\
&= \int K_{\mathcal{J}}(y-x)\left(\int K_{\mathcal{J}_0^x}(t-y)[g(t)-g(y)]\,dt\right)dy \\
&= \int K_{\mathcal{J}}(v)\left[\int K_{\mathcal{J}_0^x}(z)(g(z+v+x)-g(v+x))\,dz\right]dv \\
&= \int \mathsf{K}_{(\mathcal{A},\lambda)}(M_\vartheta^T v)\int \mathsf{K}_{(\mathcal{A},\lambda)}(M_{\vartheta^x}^T z)(g(z+v+x)-g(v+x))\,dz\,dv.
\end{aligned}
$$

Define $G_x(\cdot) = G(\cdot + x)$ and $f_x(\cdot) = f(\cdot + G(x))$. Then $g(z+v+x) = f(G_x(z+v))$ and $g(v+x) = f(G_x(v))$. Note that, for all $x \in [-1,1]^d$,

$$(63) \qquad G_x \in \mathbb{H}_d(\beta, L_2), \qquad f_x \in \mathbb{H}_1(\gamma, L_1).$$

If $1 < \gamma \leq 2$, the second property in (63) implies

$$(64) \qquad f'_x \in \mathbb{H}_1(\gamma - 1, 2L_1).$$

In the case where $1 < \beta \leq 2$, for all $u \in \mathbb{R}^d, x \in [-1,1]^d$ we define $\tilde{G}_x(u) = G_x(u) - G_x(0) - [\nabla G_x(0)]^T u$. In view of (63), for all $x \in [-1,1]^d$ we have

$$(65) \qquad \|\nabla \tilde{G}_x(u)\| \leq 2L_2 \qquad \forall u \in \mathbb{R}^d,$$

$$(66) \quad |\tilde{G}_x(t) - \tilde{G}_x(u) - [\nabla \tilde{G}_x(u)]^T(t-u)| \leq L_2\|t-u\|^\beta \qquad \forall t,u \in \mathbb{R}^d,$$

$$\Rightarrow |\tilde{G}_x(u)| \leq L_2\|u\|^\beta, \qquad u \in \mathbb{R}^d.$$

It follows from the definition of $\mathsf{K}_{(\mathcal{A},\lambda)}$ and Lemmas 1–4 that

$$(67) \qquad \int \|v\|^{\gamma\beta}|\mathsf{K}_{(\mathcal{A},\lambda)}(v)|\,dv \leq c'_6 \lambda \qquad \forall \mathcal{A} \in (0,2]^2, \lambda > 0,$$

where $c'_6 > 0$ is a constant depending only on $\mathcal{L}$ and $d$. Furthermore, for any $\mathcal{A} = (\gamma, \beta) \in (0,2]^2$ and any $\lambda \leq 1$ the support of $\mathsf{K}_{(\mathcal{A},\lambda)}$ is contained in a ball $\{u \in \mathbb{R}^d : \|u\| \leq c_K \lambda^{1/(\gamma\beta)}\}$ where the constant $c_K > 0$ depends only on $d$. Therefore,

$$(68) \quad \mathsf{K}_{(\mathcal{A},\lambda)}(M_\vartheta^T u) = 0 \qquad \forall u, \vartheta \in \mathbb{R}^d : \|u\| > c_K \lambda^{1/(\gamma\beta)}, \qquad \|\vartheta\| = 1. \qquad \square$$

2°. PROOF FOR THE ZONE OF RISD LOCAL MODEL: $1 < \gamma \leq \beta \leq 2$. Using (63) and the Taylor expansion for $G_x$ we obtain, for all $x \in [-1,1]^d$, $z, v \in \mathbb{R}^d$,

$$(69) \quad \begin{aligned} g(z+v+x) &= f(G_x(0) + [\nabla G_x(0)]^T(z+v) + \tilde{G}_x(z+v)) \\ &= f_x([\nabla G_x(0)]^T(z+v) + \tilde{G}_x(z+v)). \end{aligned}$$

Note that, by definition, $\nabla G_x(0) = \nabla G(x) = \vartheta_0^x \|\nabla G(x)\|$. Set $\nabla G_* = \vartheta^x \|\nabla G(x)\|$ and define

$$g_*(z+v+x) = f_x([\nabla G_*]^T(z+v) + \tilde{G}_x(z+v)).$$

We now approximate $g(z+v+x)$ by $g_*(z+v+x)$ in the last line of (62). In view of (68), it suffices to consider there only the values $z, v$ satisfying $\|z\|, \|v\| \leq c_K$. For such $z, v$ and all $x \in [-1, 1]^d$, the condition $\|\vartheta_0^x - \vartheta^x\| \leq \varepsilon$ and (63) imply

(70) $\quad |g(z+v+x) - g_*(z+v+x)| \leq 2c_K L_1 \|\nabla G(x)\| \varepsilon \leq 2c_K L_1 L_2 \varepsilon.$

Using (63)–(66), the Taylor expansion for $f_x$ and (64), we get that for all $x \in [-1, 1]^d$, $z, v \in \mathbb{R}^d$ the following representation holds:

(71)
$$\begin{aligned}
g_*(z+v+x) &= f_x([\nabla G_*]^T(z+v)) \\
&\quad + f'_x([\nabla G_*]^T(z+v))\tilde{G}_x(z+v) + B_{x,1}(z,v)\|z+v\|^{\gamma\beta} \\
&= f_x([\nabla G_*]^T(z+v)) \\
&\quad + [f'_x([\nabla G_*]^T(z+v)) - f'_x([\nabla G_*]^T v)] \\
&\quad \times (\tilde{G}_x(v) + [\nabla \tilde{G}_x(v)]^T z) \\
&\quad + f'_x([\nabla G_*]^T v)(\tilde{G}_x(z+v) - \tilde{G}_x(v)) \\
&\quad + f'_x([\nabla G_*]^T v)\tilde{G}_x(v) \\
&\quad + B_{x,2}(z,v)|[\nabla G_*]^T z|^{\gamma-1}\|z\|^\beta + B_{x,1}(z,v)\|z+v\|^{\gamma\beta},
\end{aligned}$$

where, for all $x \in [-1, 1]^d$, $z, v \in \mathbb{R}^d$, $B_{x,1}(\cdot, \cdot)$ and $B_{x,2}(\cdot, \cdot)$ are functions satisfying

(72) $\qquad |B_{x,1}(z,v)| \leq L_1 L_2^\gamma, \qquad |B_{x,2}(z,v)| \leq 2L_1 L_2.$

Putting $z = 0$ in (71) we obtain

(73) $\quad g_*(v+x) = f_x([\nabla G_*]^T v) + f'_x([\nabla G_*]^T v)\tilde{G}_x(v) + B_{x,1}(0,v)\|v\|^{\gamma\beta}.$

From (71) and (73) we get, for all $x \in [-1, 1]^d$, $z, v \in \mathbb{R}^d$,

(74)
$$\begin{aligned}
g_*(z&+v+x) - g_*(v+x) \\
&= f_x([\nabla G_*]^T(z+v)) - f_x([\nabla G_*]^T v) \\
&\quad + [f'_x([\nabla G_*]^T(z+v)) - f'_x([\nabla G_*]^T v)](\tilde{G}_x(v) + [\nabla \tilde{G}_x(v)]^T z) \\
&\quad + f'_x([\nabla G_*]^T v)(\tilde{G}_x(z+v) - \tilde{G}_x(v)) \\
&\quad + B_{x,2}(z,v)|[\nabla G_*]^T z|^{\gamma-1}\|z\|^\beta + B_{x,1}(z,v)\|z+v\|^{\gamma\beta} \\
&\quad - B_{x,1}(0,v)\|v\|^{\gamma\beta}.
\end{aligned}$$

Put $u = M_{\vartheta^x}^T v, s = M_{\vartheta^x}^T z$. We get from (74) that

$$\begin{aligned}
g_*(M_{\vartheta^x}s &+ M_{\vartheta^x}u + x) - g_*(M_{\vartheta^x}u + x) \\
&= (\tilde{f}_x(s_1 + u_1) - \tilde{f}_x(u_1))
\end{aligned}$$

$$\text{(75)} \quad + A_{u,x}(s_1)(\overline{G}_x(u) + [\nabla \overline{G}_x(u)]^T s)$$
$$+ f'_x(\|\nabla G(x)\|u_1)(\overline{G}_x(s+u) - \overline{G}_x(u)) + \tilde{B}_{x,2}(s,u)|s_1|^{\gamma-1}\|s\|^\beta$$
$$+ \tilde{B}_{x,1}(s,u)\|s+u\|^{\gamma\beta} - \tilde{B}_{x,1}(0,u)\|u\|^{\gamma\beta},$$

where $s_1$ and $u_1$ are the first components of $s \in \mathbb{R}^d$ and $u \in \mathbb{R}^d$, respectively,

$$\tilde{f}_x(u_1) = f_x(\|\nabla G(x)\|u_1), \qquad \overline{G}_x(u) = \tilde{G}_x(M_{\vartheta^x}u),$$
$$\tilde{B}_{x,1}(s,u) = B_{x,1}(M_{\vartheta^x}s, M_{\vartheta^x}u)$$
$$\tilde{B}_{x,2}(s,u) = \|\nabla G(x)\|^{\gamma-1} B_{x,2}(M_{\vartheta^x}s, M_{\vartheta^x}u)$$

and

$$A_{u,x}(s_1) = f'_x(\|\nabla G(x)\|(s_1+u_1)) - f'_x(\|\nabla G(x)\|u_1).$$

It is easy to see that inequalities (65) and (66) remain valid with $\overline{G}_x$ in place of $\tilde{G}_x$.

Now for all $x \in [-1,1]^d$, $s, u \in \mathbb{R}^d$ we introduce

$$q_{u,x}(s_1) = (\tilde{f}_x(s_1+u_1) - \tilde{f}_x(u_1)) + A_{u,x}(s_1)(\overline{G}_x(u) + [\nabla \overline{G}_x(u)]^T \vartheta^x s_1)$$
$$+ f'_x(\|\nabla G(x)\|u_1)[\nabla \overline{G}_x(u)]^T \vartheta^x s_1,$$
$$p_{u,x}(s) = f'_x(\|\nabla G(x)\|u_1)(\overline{G}_x(s+u) - \overline{G}_x(u) - [\nabla \overline{G}_x(u)]^T s),$$
$$B^{u,x}(s) = \tilde{B}_{x,2}(s,u),$$
$$Q_{u,x}(s) = q_{u,x}(s_1) + p_{u,x}(s) + \tilde{B}_{x,2}(s,u)|s_1|^{\gamma-1}\|s\|^\beta,$$
$$P_{u,x}(s) = f'_x(\|\nabla G(x)\|(s_1+u_1))[\nabla \overline{G}_x(u)]^T s_\perp,$$

where $s_\perp = s - s_1 \vartheta^x$. With this notation (75) can be written as

$$\text{(76)} \quad \begin{aligned} g_*(M_{\vartheta^x}s + M_{\vartheta^x}u + x) &- g_*(M_{\vartheta^x}u + x) \\ &= Q_{u,x}(s) + P_{u,x}(s) + \tilde{B}_{x,1}(s,u)\|s+u\|^{\gamma\beta} - \tilde{B}_{x,1}(0,u)\|u\|^{\gamma\beta}. \end{aligned}$$

We now prove that, for all $x \in [-1,1]^d$ and all $u \in \mathbb{R}^d$ such that $\|u\| \leq c_K \lambda^{1/(\gamma\beta)}$ [cf. (68)], the triplet $(q_{u,x}, p_{u,x}, B^{u,x})$ belongs to the set $\mathfrak{B}(\mathcal{A}, \lambda)$ (cf. definition before Lemma 3), and thus Lemmas 3 or 4 can be applied. We need to check (23)–(25).

*Checking* (23). In view of (63) we have

$$|\tilde{f}_x(s_1+u_1) - \tilde{f}_x(u_1) - \tilde{f}'_x(u_1)s_1| \leq L_1 L_2 |s_1|^\gamma.$$

Therefore,

$$\text{(77)} \quad \left| \frac{1}{2\lambda_0^{1/\gamma}} \int_{-\lambda_0^{1/\gamma}}^{\lambda_0^{1/\gamma}} (\tilde{f}_x(s_1+u_1) - \tilde{f}_x(u_1))\, ds_1 \right| \leq \frac{L_1 L_2}{2\lambda_0^{1/\gamma}} \int_{-\lambda_0^{1/\gamma}}^{\lambda_0^{1/\gamma}} |s_1|^\gamma \, ds_1$$
$$\leq \frac{L_1 L_2}{2} \lambda.$$

Next, remark that (64) implies $|A_{u,x}(s_1)| \leq 2L_1 L_2^{\gamma-1}|s_1|^{\gamma-1}$. Furthermore, (66) with $\overline{G}_x$ in place of $\tilde{G}_x$ yields $|\overline{G}_x(u)| \leq L_2\|u\|^\beta$. Now, $q_{u,x}(0) = 0$ and using these remarks, (77) and (65) we get, for $\|u\| \leq c_K \lambda^{1/(\gamma\beta)}$,

$$
\begin{aligned}
\left|\frac{1}{2\lambda_0^{1/\gamma}} \int_{-\lambda_0^{1/\gamma}}^{\lambda_0^{1/\gamma}} q_{u,x}(s_1)\, ds_1\right| \\
\leq \frac{L_1 L_2}{2}\lambda + \frac{1}{2\lambda_0^{1/\gamma}} \int_{-\lambda_0^{1/\gamma}}^{\lambda_0^{1/\gamma}} |A_{u,x}(s_1)|(|\overline{G}_x(u)| + \|\nabla \overline{G}_x(u)\|\|s_1\|)\, ds_1 \\
\leq \frac{L_1 L_2}{2}\lambda + 2 L_1 L_2^\gamma \left(\frac{1}{\gamma}\lambda^{(\gamma-1)/\gamma}\|u\|^\beta + \frac{2}{\gamma+1}\lambda\right) \\
\leq \left[\frac{L_1 L_2}{2} + 2 L_1 L_2^\gamma \left(\frac{(2c_K)^\beta}{\gamma} + \frac{2}{\gamma+1}\right)\right]\lambda \leq c_3 \lambda,
\end{aligned}
$$
(78)

where the constant $c_3$ depends only on $\mathcal{L}$ and $d$. It can be taken as a maximum of the last expression in square brackets over $(\gamma, \beta) \in [1,2]^2$.

*Checking* (24) *and* (25). It suffices to note that, for all $x \in [-1,1]^d$, the first property in (66) with $\overline{G}_x$ in place of $\tilde{G}_x$ and the second property in (63) yield

$$|p_{u,x}(s') - p_{u,x}(s) - [\nabla p_{u,x}(s)]^T(s'-s)| \leq |f'_x(\|\nabla G(x)\|u_1)|L_2\|s'-s\|^\beta$$
$$\leq L_1 L_2 \|s'-s\|^\beta \quad \forall s, s' \in \mathbb{R}^d.$$

This proves (24) with $b = \beta$ and $L = L_1 L_2$. Finally, (25) with $B = B^{u,x}$, $c_4 = 2L_1 L_2^\gamma$ follows from (72).

We are now in a position to apply Lemmas 3 and 4. We demonstrate this, for example, for Lemma 4. Take there $q = q_{u,x}, p = p_{u,x}, B = B^{u,x}$ for any $\|u\| \leq c_K \lambda^{1/(\gamma\beta)}$ and $x \in [-1,1]^d$. Since $Q_{u,x}(0) = 0$, the result (29) of Lemma 4 yields

$$\left|\int \mathsf{K}_{(\mathcal{A},\lambda)}(s) Q_{u,x}(s)\, ds\right| \leq c_5 \lambda, \tag{79}$$

where $c_5$ depends only on $\mathcal{L}$ and $d$. Furthermore, by construction the weight $\mathsf{K}_{(\mathcal{A},\lambda)}$ is symmetric, that is, $\mathsf{K}_{(\mathcal{A},\lambda)}(s) = \mathsf{K}_{(\mathcal{A},\lambda)}(-s)$ and hence

$$\int \mathsf{K}_{(\mathcal{A},\lambda)}(s) P_{u,x}(s)\, ds = 0. \tag{80}$$

Next, using (72) we find

$$|\tilde{B}_{x,1}(s,u)\|s+u\|^{\gamma\beta} - \tilde{B}_{x,1}(0,u)\|u\|^{\gamma\beta}| \leq 2^{\gamma\beta} L_1 L_2^\gamma (\|s\|^{\gamma\beta} + \|u\|^{\gamma\beta}).$$

Combining this inequality and (79)–(80) with (76) we get, for all $x \in [-1,1]^d$, $u \in \mathbb{R}^d$,

$$\left| \int \mathsf{K}_{(\mathcal{A},\lambda)}(s)(g_*(M_{\vartheta^x}s + M_{\vartheta^x}u + x) - g_*(M_{\vartheta^x}u + x))\, ds \right|$$

$$\leq c_5\lambda + 2^{\gamma\beta} L_1 L_2^\gamma \left[ \int |\mathsf{K}_{(\mathcal{A},\lambda)}(s)| \|s\|^{\gamma\beta}\, ds + \|\mathsf{K}_{(\mathcal{A},\lambda)}\|_1 \|u\|^{\gamma\beta} \right].$$

We finally get (21) from this inequality invoking (67), (62), (70) and recalling that $\|\mathsf{K}_{(\mathcal{A},\lambda)}\|_1 = \|K_{\mathcal{J}}\|_1$ for all $\mathcal{A} \in (0,2]^2, \lambda > 0$, and $\|\mathsf{K}_{(\mathcal{A},\lambda)}\|_1 = \|K_{\mathcal{J}_0^x}\|_1$. □

$3°$. PROOF OF (21) FOR THE LOCAL SINGLE-INDEX ZONE: $\gamma \leq 1, 1 < \beta \leq 2$. Using (66) and the second property in (63), for all $z, v \in \mathbb{R}^d$, $x \in [-1,1]^d$ we may write

$$g_*(z + v + x) = f_x([\nabla G_*]^T(z+v)) + B_{x,1}(z,v)\|z+v\|^{\gamma\beta},$$

where $B_{x,1}$ satisfies (72). This can be viewed as a simplified version of (71). Following almost the same argument as in $2°$ (the main difference is that now we drop all the terms containing $f'_x$ and $B_{x,2}$) and applying Lemma 2 we obtain (21). □

$4°$. PROOF OF (21) FOR THE ZONE OF SLOW RATE: $(\gamma, \beta) \in (0,1]^2$. Using the Hölder condition on $f$ and $G_x$ we obtain, for all $z, v \in \mathbb{R}^d, x \in [-1,1]^d$,

$$g(z + v + x) \equiv f(G_x(z+v)) = f(G_x(0)) + B_{x,1}(z,v)\|z+v\|^{\gamma\beta},$$

where $B_{x,1}$ satisfies (72). Now, (21) easily follows from this relation, (62), (67) and the definition of $\mathsf{K}_{(\mathcal{A},\lambda)}$ for the zone of slow rate. □

$5°$. PROOF OF (21) FOR THE ZONE OF INACTIVE STRUCTURE: $1 < \beta \leq \gamma \leq 2$. Since $f \in \mathbb{H}_1(\gamma, L_1)$ and $\|\nabla G_x(\cdot)\| \leq L_2$, for all $z, v \in \mathbb{R}^d, x \in [-1,1]^d$ we may write

$$f(G_x(z+v)) = f(G_x(v)) + f'(G_x(v))(G_x(z+v) - G_x(v)) + B_{x,1}(z,v)\|z\|^\gamma$$

$$= f(G_x(v)) + f'(G_x(v))(G_x(z+v) - G_x(v) - [\nabla G_x(v)]^T z)$$

$$+ f'(G_x(v))[\nabla G_x(v)]^T z + B_{x,1}(z,v)\|z\|^\gamma$$

$$= f(G_x(v)) + f'(G_x(v))[\nabla G_x(v)]^T z + B_{x,2}(z,v)\|z\|^\beta$$

$$+ B_{x,1}(z,v)\|z\|^\gamma,$$

where $B_{x,1}$ satisfies (72) and $|B_{x,2}(\cdot,\cdot)| \leq L_1 L_2$. Since the weight $\mathsf{K}_{(\mathcal{A},\lambda)}$ is symmetric,

$$\int \mathsf{K}_{(\mathcal{A},\lambda)}(M_{\vartheta^x}^T z) f'(G_x(v))[\nabla G_x(v)]^T z\, dz = 0.$$

Now, (21) easily follows from these relations, (62), the definition of $\mathsf{K}_{(\mathcal{A},\lambda)}$ for the zone of inactive structure and the condition $\lambda \leq 1$. □

6°. PROOF OF (22). For a function $K \in L_2(\mathbb{R}^d)$, let us denote by $\widehat{K}$ its Fourier transform. Using Parceval's identity we obtain, for any $\mathcal{J}, \mathcal{J}' \in \mathfrak{J}$,

$$\|\Delta_{\mathcal{J}'} K_{\mathcal{J}' * \mathcal{J}}\|_2 = \frac{1}{\sqrt{2\pi}} \|\widehat{\Delta_{\mathcal{J}'} K}_{\mathcal{J}' * \mathcal{J}}\|_2 = \frac{1}{\sqrt{2\pi}} \|(\widehat{K}_{\mathcal{J}} - 1)\widehat{K}_{\mathcal{J}'}\|_2$$

$$\leq \frac{1}{\sqrt{2\pi}} (\|\widehat{K}_{\mathcal{J}}\|_\infty + 1)\|\widehat{K}_{\mathcal{J}'}\|_2 \leq (\|K_{\mathcal{J}}\|_1 + 1)\|K_{\mathcal{J}'}\|_2.$$

Since $\int K_{\mathcal{J}'} = 1$, this proves (22). □

**A.2. Proof of Lemma 3.** First, note that some cases are trivial because the number $r$ of steps of the weight construction is bounded by 3. In fact, if $(\rho+1)\rho < (\beta-\gamma)/\gamma$ and $V(\lambda) \leq \ln(\frac{\sqrt{5}+1}{2})$ we have $r \leq 3$ by definition. If $(\rho+1)\rho \geq (\beta-\gamma)/\gamma$ we use the weight as in Lemma 3. But for this weight the condition $(\rho+1)\rho \geq (\beta-\gamma)/\gamma$ implies that, again, $r \leq 3$.

So, we will treat only the remaining case where $(\rho+1)\rho < (\beta-\gamma)/\gamma$ and $V(\lambda) > \ln(\frac{\sqrt{5}+1}{2})$. The last inequality implies that $r > 3$.

Note that, by definition, $\alpha < \frac{1}{2}\ln(\frac{\sqrt{5}+1}{2})$. Further, for $r \geq 3$ we have also the lower bound: $\alpha \geq \frac{1}{4}\ln(\frac{\sqrt{5}+1}{2})$. Thus for $r \geq 3$,

$$(81) \qquad 0.786 \leq \left(\frac{\sqrt{5}+1}{2}\right)^{-1/2} < e^{-\alpha} \leq \left(\frac{\sqrt{5}+1}{2}\right)^{-1/4} \leq 0.887.$$

1°. PROOF OF (29). From the definition of $\mathsf{K}_{(\mathcal{A},\lambda)}$ we find

$$[\mathsf{K}_{(\mathcal{A},\lambda)} * q](0) = 2^{-d} \sum_{i=1}^{r} \int \Lambda_i(|y|)q(y_1)\,dy = 2^{-d}\int \Lambda_1(|y|)q(y_1)\,dy$$

$$= \frac{1}{u_1}\int \frac{q(y_1)+q(-y_1)}{2}\mathbb{I}_{[0,u_1]}(y_1)\,dy_1,$$

where $u_1 = \lambda^{1/\gamma}$. This and (23) imply

$$(82) \quad |[\mathsf{K}_{(\mathcal{A},\lambda)} * q](0) - q(0)| = \left|(2\lambda^{1/\gamma})^{-1}\int_{-\lambda^{1/\gamma}}^{\lambda^{1/\gamma}} q(y_1)\,dy_1 - q(0)\right| \leq c_3\lambda.$$

We now obtain a similar bound for $|[\mathsf{K}_{(\mathcal{A},\lambda)} * p](0) - p(0)|$. Note that, in view of (24), for all $z = (z_1, \ldots, z_d) \in \mathbb{R}^d$ we have

$$(83) \qquad p(z) = \tilde{p}(z) + z_1 \frac{\partial p}{\partial z_1}(0, z_2, \ldots, z_d) + B_1(z)z_1^\beta,$$

where $\tilde{p}(z) = p(0, z_2, \ldots, z_d)$ and $\sup_{z \in \mathbb{R}^d} |B_1(z)| \leq L$. For the same reason, for all $z_{(d-1)} \triangleq (0, z_2, \ldots, z_d)$ we have

$$\tilde{p}(z) = \tilde{p}(0) + [\nabla \tilde{p}(0)]^T z_{(d-1)} + B_2(z_{(d-1)}) \|z_{(d-1)}\|^\beta, \tag{84}$$

where as previously $|B_2(\cdot)| \leq L$. Combining (83) and (84) and taking into account that the function $\mathsf{K}_{(\mathcal{A},\lambda)}$ is symmetric, $\int \mathsf{K}_{(\mathcal{A},\lambda)} = 1$ and $\tilde{p}(0) = p(0)$ we get

$$\begin{aligned}
&|[\mathsf{K}_{(\mathcal{A},\lambda)} * p](0) - p(0)| \\
&= \left| \int \mathsf{K}_{(\mathcal{A},\lambda)}(z)(B_1(z)z_1^\beta + B_2(z_{(d-1)})\|z_{(d-1)}\|^\beta) \, dz \right|.
\end{aligned} \tag{85}$$

Now

$$\begin{aligned}
&\left| \int \mathsf{K}_{(\mathcal{A},\lambda)}(z) B_2(z_{(d-1)}) \|z_{(d-1)}\|^\beta \, dz \right| \\
&= \left| (2(v_1 - v_2))^{1-d} \int B_2(z_{(d-1)}) \|z_{(d-1)}\|^\beta \mathbb{I}_{[v_2,v_1]^{d-1}}(|z_{(d-1)}|) \, dz_{(d-1)} \right. \\
&\quad + \sum_{i=1}^{r-1} \Big[ (2(v_i - v_{i+1}))^{1-d} \\
&\qquad \times \int B_2(z_{(d-1)}) \|z_{(d-1)}\|^\beta \mathbb{I}_{[v_{i+1},v_i]^{d-1}}(|z_{(d-1)}|) \, dz_{(d-1)} \\
&\qquad - (2(v_{i-1} - v_i))^{1-d} \int B_2(z_{(d-1)}) \|z_{(d-1)}\|^\beta \\
&\qquad\qquad \left. \times \mathbb{I}_{[v_i,v_{i-1}]^{d-1}}(|z_{(d-1)}|) \, dz_{(d-1)} \Big] \right| \\
&\leq (2v_r)^{1-d} \int |B_2(z_{(d-1)})| \|z_{(d-1)}\|^\beta \mathbb{I}_{[0,v_r]^{d-1}}(|z_{(d-1)}|) \, dz_{(d-1)} \\
&= (\lambda^{1/\beta})^{1-d} \int |B_2(z_{(d-1)})| \|z_{(d-1)}\|^\beta \mathbb{I}_{[0,\lambda^{1/\beta}]^{d-1}}(|z_{(d-1)}|) \, dz_{(d-1)} \\
&\leq 2^{d-1} d^{\beta/2} L\lambda \leq 2^{d-1} dL\lambda,
\end{aligned} \tag{86}$$

where $|z_{(d-1)}| = (|z_2|, \ldots, |z_d|)$. Further, note that $v \geq u \geq 1$ implies $e^{v/u} \leq e^v/u$ [in fact, $v(1 - 1/u) \geq u - 1 \geq \ln u$]. Using this remark and the fact that $\frac{\beta}{\gamma - 1} > 1$ we find

$$\begin{aligned}
u_i &= \lambda^{1/\gamma} \exp\left( \frac{\beta}{\gamma - 1} \exp(\alpha(i-1)) \right) = \lambda^{1/\gamma} \exp\left( \frac{\beta}{\gamma - 1} \exp(\alpha i) e^{-\alpha} \right) \\
&\leq u_{i+1} e^{-\alpha}, \qquad i = 1, \ldots, r-1
\end{aligned} \tag{87}$$

and therefore $u_i/u_r \leq e^{\alpha(i-r)}$. This and the equality $u_r = \lambda^{1/\beta}$ allow us to get

$$
\begin{aligned}
(88) \quad &\left|\int \mathsf{K}_{(\mathcal{A},\lambda)}(z) B_1(z) z_1^\beta \, dz\right| \\
&\leq L \int |\mathsf{K}_{(\mathcal{A},\lambda)}(z)| |z_1|^\beta \, dz \\
&= \frac{L}{u_1} \int z_1^\beta \mathbb{I}_{[0,u_1]}(z_1) \, dz_1 + \sum_{i=2}^{r} \frac{2L}{u_i - u_{i-1}} \int z_1^\beta \mathbb{I}_{[u_{i-1},u_i]}(z_1) \, dz_1 \\
&\leq 2L \sum_{i=1}^{r} u_i^\beta \leq 2L\lambda \sum_{i=1}^{r} \left(\frac{u_i}{u_r}\right)^\beta \leq 2\lambda L \sum_{l=0}^{\infty} e^{-\alpha l} = 2\lambda L (1 - e^{-\alpha})^{-1}.
\end{aligned}
$$

From (85), (86) and (88) we get

$$
(89) \quad |[\mathsf{K}_{(\mathcal{A},\lambda)} * p](0) - p(0)| \leq \lambda L [2^{d-1} d + 2(1 - e^{-\alpha})^{-1}].
$$

We now estimate the value $|\int \mathsf{K}_{(\mathcal{A},\lambda)}(y) B(y) y_1^{\gamma-1} \|y\|^\beta \, dy|$. In view of (42),

$$
\begin{aligned}
(90) \quad &u_1^{\gamma-1} v_1^\beta \leq \lambda \exp\{\beta - \nu\beta e^\alpha\} \leq \lambda \exp\{(1-\nu)\beta\}, \\
&u_i^{\gamma-1} v_i^\beta \leq u_i^{\gamma-1} v_{i-1}^\beta = \lambda \exp\{(1-\nu)\beta \exp(\alpha(i-1))\}, \qquad i = 2, \ldots, r.
\end{aligned}
$$

Using (90), we get similarly to (88):

$$
\begin{aligned}
(91) \quad &\left|\int \mathsf{K}_{(\mathcal{A},\lambda)}(y) B(y) y_1^{\gamma-1} \|y\|^\beta \, dy\right| \\
&\leq c_4 \int |\mathsf{K}_{(\mathcal{A},\lambda)}(y)| |y_1|^{\gamma-1} \sum_{j=1}^{d} |y_j|^\beta \, dy \\
&= c_4 \left[\int |\mathsf{K}_{(\mathcal{A},\lambda)}(y)| |y_1|^{\gamma+\beta-1} \, dy + \sum_{j=2}^{d} \int |\mathsf{K}_{(\mathcal{A},\lambda)}(y)| |y_1|^{\gamma-1} |y_j|^\beta \, dy\right] \\
&\leq 2c_4 \left[\sum_{i=1}^{r} u_i^{\beta+\gamma-1} + d \sum_{i=1}^{r} u_i^{\gamma-1} v_i^\beta\right] \\
&\leq 2c_4 \left[\lambda^{(\beta+\gamma-1)/\beta} \sum_{l=0}^{\infty} e^{-\alpha l (\beta+\gamma-1)} + \lambda d \sum_{l=0}^{\infty} \exp\{(1-\nu)\beta \exp(\alpha l)\}\right] \\
&\leq 2c_4 \lambda [(1 - e^{-\alpha})^{-1} + d(1 - e^{(1-\nu)\alpha})^{-1}],
\end{aligned}
$$

where the last inequality holds for $0 < \lambda \leq 1$ and we used that $\beta \exp(\alpha l) \geq \alpha l$, $\nu > 1$. Summing up the results of (82), (89), (91) and taking into account (81) we obtain (29). $\square$

2°. PROOF OF (30). In the same way as above we get, for $0 < \lambda \leq 1$,

$$\int |\mathsf{K}_{(\mathcal{A},\lambda)}(y)| \|y\|^m \, du \leq d^{m/2} \int |\mathsf{K}_{(\mathcal{A},\lambda)}(y)| \sum_{j=1}^d |y_j|^m \, dy$$

$$\leq 2d^{m/2} \left[ \sum_{i=1}^r u_i^m + d \sum_{i=1}^r v_i^m \right]$$

$$\leq C(d) \lambda^{m/(\gamma\beta)} [(1 - e^{-m\alpha})^{-1} + (1 - e^{m\nu\alpha})^{-1}].$$

Here and in what follows use the same notation $C(d)$ for possibly different positive constants depending only on $d$. □

3°. PROOF OF (31). Since $\nu < 2 < \frac{\beta}{\beta-\gamma}$ we have, for $0 < \lambda \leq 1$,

$$v_{r-1} \triangleq \lambda^{1/(\gamma\beta)} \exp\{-\nu \exp(\alpha(r-1))\} = \lambda^{1/(\gamma\beta) + \nu(\gamma-1)(\beta-\gamma)/(\gamma\beta)^2} \geq \lambda^{1/\beta}.$$

By the definition of $v_r$ this implies that $v_{r-1} - v_r \geq \lambda^{1/\beta}/2$. Further, as $u_r = \lambda^{1/\beta}$, in view of (87), we have

$$u_r - u_{r-1} \geq (1 - e^{-\alpha})\lambda^{1/\beta}.$$

We deduce that

(92) $$\mu_{r,r-1} \geq \mu_{r,r} \geq 2^{1-d} \lambda^{d/\beta}(1 - e^{-\alpha}).$$

Note that by (87),

$$u_{i+1} - u_i \geq (1 - e^{-\alpha})u_{i+1} \qquad \text{for } i = 1, \ldots, r-1.$$

Also, as $\nu > 1$, it is straightforward to check that

$$v_i - v_{i+1} \geq (1 - e^{-\alpha})v_i \qquad \text{for } i = 1, \ldots, r-2.$$

Thus, we get

(93) $$\mu_{1,1} = u_1(v_1 - v_2)^{d-1} \geq (1 - e^{-\alpha})^{d-1} \exp(-(d-1)\nu e^\alpha)\lambda^{1/\gamma + (d-1)/\beta}.$$

Recall that we are considering the case where $\rho(1+\rho) < (\beta-\gamma)/\gamma, 1 < \gamma \leq \beta \leq 2$, so that $\rho(1+\rho) < 1$, and thus $\rho < \frac{\sqrt{5}-1}{2}$. This and the choice of parameters $\alpha, \nu$ combined with (81) implies

$$e^{-\alpha} - \rho\nu \geq \left(\frac{\sqrt{5}+1}{2}\right)^{-1/2} - \rho\nu \geq \left(\frac{\sqrt{5}+1}{2}\right)^{-1/2} - \frac{\sqrt{5}-1}{2}\nu \triangleq \delta \geq 0.0891.$$

Now,

$$\frac{\beta}{\gamma - 1} e^{-\alpha} - (d-1)\nu \geq \frac{\delta\beta}{\gamma - 1} \geq 2\delta.$$

Hence, for $i = 2, \ldots, r-1$ we have

$$\mu_{i,i-1} \geq \mu_{i,i}$$
$$\geq C(d)\lambda^{1/\gamma + (d-1)/(\gamma\beta)}$$
(94)
$$\times \exp\left\{\frac{\beta}{\gamma - 1}\exp(\alpha(i-1)) - (d-1)\nu\exp(\alpha i)\right\}$$
$$\geq C(d)\lambda^{1/\gamma + (d-1)/(\gamma\beta)}\exp\{2\delta\exp(\alpha i)\}.$$

Note that

$$(95) \quad \|\mathsf{K}_{(\mathcal{A},\lambda)}\|_2^2 = \mu_{1,1}^{-1} + \sum_{i=2}^r (\mu_{i,i-1}^{-1} + \mu_{i,i}^{-1}) \leq \mu_{1,1}^{-1} + 2\sum_{i=2}^r \mu_{i,i}^{-1}.$$

We deduce from (92)–(95) that

$$\|\mathsf{K}_{(\mathcal{A},\lambda)}\|_2^2 \leq C(d)(\lambda^{1/\gamma + (d-1)/(\gamma\beta)} + \lambda^{-d/\beta}).$$

This proves the second inequality in (31). The first inequality becomes obvious if we note that $V(\lambda) \leq \ln\ln(1/\lambda)$ and so $\|\mathsf{K}_{(\mathcal{A},\lambda)}\|_1 = 2r - 1 \leq c_7 \ln\ln(1/\lambda)$, for $\lambda$ small enough, where $c_7$ is an absolute constant. $\square$

**A.3. Proof of Lemma 3.** Following the same lines as in the proof of (29) in Lemma 4 we obtain the bound (26) of Lemma 3 with

$$c_5 = C(d)(c_3 + Lr + c_4 r).$$

1°. PROOF OF (27). By definition, $u_r = \lambda^{1/\beta}$ and for $0 < \lambda \leq 1$ we have $u_2 \geq \lambda^{1/\gamma}$, so that $v_1 = \lambda^{1/\beta} u_2^{-(\gamma-1)/\beta} \leq \lambda^{1/(\gamma\beta)}$. Using these remarks and acting as in the proof of (30) in Lemma 4 we obtain, for $0 < \lambda \leq 1$,

$$\int |\mathsf{K}_{(\mathcal{A},\lambda)}(y)|\|y\|^m\,du \leq 2d^{m/2}\left[\sum_{i=1}^r u_i^m + d\sum_{i=1}^r v_i^m\right]$$
$$\leq 2d^{m/2}r(u_r^m + dv_1^m) \leq C(d)r\lambda^{m/(\gamma\beta)}. \quad \square$$

2°. PROOF OF (28). Observe that $\alpha_{j+1} - \alpha_j > 0$ for $j = 1, \ldots, r-1$, so that for $\lambda \to 0$ we have $u_j/u_{j-1} \to \infty$ and $v_{j-1}/v_j \to \infty$. In particular,

$$\mu_{j,j-1} = (u_j - u_{j-1})(v_{j-1} - v_j)^{d-1} \geq \mu_{j,j} = (u_j - u_{j-1})(v_j - v_{j+1})^{d-1}$$
$$\geq \tfrac{1}{2}u_j v_j^{d-1}$$

for all $\lambda$ small enough. Next note that, by definition,

$$\alpha_{r-2} \geq (\alpha_{r-1} - \beta^{-1})\rho^{-1} \geq \frac{\beta - \gamma}{\gamma\beta\rho}.$$

Then $u_2 \leq \lambda^{(\beta-\gamma)/(\gamma\beta\rho)}$ and for $\lambda$ small enough we get by the definition of $\rho$:

$$\mu_{1,1} \geq \tfrac{1}{2} u_1 v_1^{d-1} = \tfrac{1}{2}\lambda^{(d-1)/\beta} u_1 u_2^{-\rho} = \tfrac{1}{2}\lambda^{(d-1)/\beta}\lambda^{1/\gamma - (\beta-\gamma)/(\gamma\beta)} = \tfrac{1}{2}\lambda^{d/\beta}.$$

Further, as $u_r = \lambda^{1/\beta}$ and $v_r = \tfrac{1}{2}\lambda^{1/\beta}$, $v_{r+1} = 0$,

$$\mu_{r,r} \geq 2^{-d}\lambda^{d/\beta}$$

for $\lambda$ small enough. Next, for $1 < j < r$,

$$\mu_{j,j} \geq \tfrac{1}{2} u_j v_j^{d-1} = \tfrac{1}{r2}\lambda^{(d-1)/\beta} u_j u_{j+1}^{-\rho}.$$

By the definition of the sequence $(\alpha_k)$,

$$(d-1)/\beta + \alpha_k - \rho/\alpha_{k-1} = d/\beta, \qquad k = 1, \ldots, r-1.$$

Thus

$$\mu_{j,j} \geq \tfrac{1}{2}\lambda^{(d-1)/\beta + \alpha_{r-j} - \rho\alpha_{r-(j+1)}} = \tfrac{1}{2}\lambda^{d/\beta}, \qquad j = 2, \ldots, r-1.$$

Substitution of the above bounds into (95) yields

$$\|\mathsf{K}_{(\mathcal{A},\lambda)}\|_2^2 \leq C(d)\lambda^{-d/\beta} r. \qquad \square$$

A. B. JUDITSKY  
LABORATOIRE JEAN KUNTZMANN  
UNIVERSITÉ GRENOBLE 1  
B.P. 53, 38041 GRENOBLE  
FRANCE  
E-MAIL: anatoli.juditsky@imag.fr  

O. V. LEPSKI  
LABORATOIRE D'ANALYSE, TOPOLOGIE ET PROBABILITÉS  
UNIVERSITÉ DE PROVENCE  
39, RUE F. JOLIOT CURIE, 3453 MARSEILLE  
FRANCE  
E-MAIL: lepski@cmi.univ-mrs.fr

A. B. Tsybakov  
Laboratoire de Statistique  
CREST  
Timbre J340  
3, av. Pierre Larousse, 92240 Malakoff  
France  
and  
Laboratoire de Probabilités et Modèles Aléatoires  
Université Paris 6  
4, pl. Jussieu  
Case 188, 75252 Paris  
France  
E-mail: alexandre.tsybakov@ensae.fr  
   tsybakov@ccr.jussieu.fr